\documentclass[11pt,a4paper]{article}

\usepackage{amsmath,amsthm,amssymb}
	
	\usepackage[utf8]{inputenc}
	\usepackage[T1]{fontenc}
	\usepackage{lmodern}
	\usepackage{mathdots}
	\usepackage{mathtools}
\usepackage{tikz}

\usepackage{color}
\usepackage{fullpage}

	\newtheorem{theorem}{Theorem}[section]
	\newtheorem{lemma}[theorem]{Lemma}
	\newtheorem{corollary}[theorem]{Corollary}
	
	\newtheorem{proposition}[theorem]{Proposition}
    \newtheorem{definition}[theorem]{Definition}
	\theoremstyle{definition}
    \DeclareMathOperator{\rank}{rank}
	\DeclareMathOperator{\spann}{span}

	\newcommand{\F}{\mathbb{F}}
	\newcommand{\C}{\mathbb{C}}

\newcommand {\mat}[1]{\left[\begin{array}{#1}}
\newcommand {\rix}          {\end{array}\right]}

	\newtheorem{remark}[theorem]{Remark}
	\usepackage{colonequals}
	\usepackage{microtype}
	\usepackage{hyperref}
	
	\DeclareMathOperator\Rev{Rev}
	
	\DeclareMathOperator\diag{diag}

\usepackage{xcolor}

	\DeclareMathOperator{\cwRev}{cwRev}

	\let\tilde\widetilde

\begin{document}

%\begin{frontmatter}

%% Title, authors and addresses

%% use the tnoteref command within \title for footnotes;
%% use the tnotetext command for the associated footnote;
%% use the fnref command within \author or \address for footnotes;
%% use the fntext command for the associated footnote;
%% use the corref command within \author for corresponding author footnotes;
%% use the cortext command for the associated footnote;
%% use the ead command for the email address,
%% and the form \ead[url] for the home page:
%%
%% \title{Title\tnoteref{label1}}
%% \tnotetext[label1]{}
%% \author{Name\corref{cor1}\fnref{label2}}
%% \ead{email address}
%% \ead[url]{home page}
%% \fntext[label2]{}
%% \cortext[cor1]{}
%% \address{Address\fnref{label3}}
%% \fntext[label3]{}

%\dochead{Root polynomials}
%% Use \dochead if there is an article header, e.g. \dochead{Short communication}

\title{The $\mathbb{DL}(P)$ vector space of pencils for singular matrix polynomials}

\author{Froil\'{a}n M. Dopico \footnote{Departamento de Matem\'aticas, Universidad Carlos III de Madrid, Avda. Universidad 30, 28911 Legan\'{e}s, Spain, \texttt{dopico@math.uc3m.es}. Supported by grant PID2019-106362GB-I00 funded by MCIN/AEI/ 10.13039/501100011033 and by the Madrid Government (Comunidad de Madrid-Spain) under the Multiannual Agreement with UC3M in the line of Excellence of University Professors (EPUC3M23), and in the context of the V PRICIT (Regional Programme of Research and Technological Innovation).} \,  and Vanni Noferini\footnote{Department of Mathematics and Systems Analysis, Aalto University, PO Box 11100, FI-00076 Aalto,
Finland. Supported by an Academy of Finland grant (Suomen Akatemian p\"{a}\"{a}tos 331240).}}

\maketitle

\begin{abstract}
Given a possibly singular matrix polynomial $P(z)$, we study how the eigenvalues, eigenvectors, root polynomials, minimal indices, and minimal bases of the pencils in the vector space $\mathbb{DL}(P)$ introduced in Mackey, Mackey, Mehl, and Mehrmann [SIAM J. Matrix Anal. Appl. 28(4), 971-1004, 2006] are related to those of $P(z)$. If $P(z)$ is regular, it is known that those pencils in $\mathbb{DL}(P)$ satisfying the generic assumptions in the so-called eigenvalue exclusion theorem are strong linearizations for $P(z)$. This property and the block-symmetric structure of the pencils in $\mathbb{DL}(P)$ have made these linearizations among the most influential for the theoretical and numerical treatment of structured regular matrix polynomials. However, it is also known that, if $P(z)$ is singular, then none of the pencils in $\mathbb{DL}(P)$ is a linearization for $P(z)$. In this paper, we prove that despite this fact a generalization of the eigenvalue exclusion theorem holds for any singular matrix polynomial $P(z)$ and that such a generalization allows us to recover all the relevant quantities of $P(z)$ from any pencil in $\mathbb{DL}(P)$ satisfying the eigenvalue exclusion hypothesis. Our proof of this general theorem relies heavily on the representation of the pencils in $\mathbb{DL} (P)$ via B\'{e}zoutians by Nakatsukasa, Noferini and Townsend [SIAM J. Matrix Anal. Appl. 38(1), 181-209, 2015].
\end{abstract}

\textbf{keyword}
B\'{e}zout matrix, \ B\'{e}zoutian, \ $\mathbb{DL}(P)$, \ singular matrix polynomial, \ minimal basis, \ minimal indices, \ root polynomial, \ eigenvalue exclusion theorem, \ linearization
%% MSC codes here, in the form: \MSC code \sep code
%% or \MSC[2008] code \sep code (2000 is the default)
%\MSC

%\end{frontmatter}
	\section{Introduction} Given an algebraically closed field $\F$ with characteristic\footnote{The assumption that the base field is closed is for simplicity of exposition, as it ensures for instance that finite eigenvalues all lie in $\F$. The assumption that $\mathrm{char}(\F)=0$ is only needed in two technical results, Lemma \ref{lem:combinatorial} and Theorem \ref{blockevaluation}, that simplify the proof of our main results. The main results of this paper could also be proved without these assumptions, but this would complicate the exposition.} $0$, this paper deals with a matrix polynomial
\begin{equation}\label{eq.intropoly}
P(z) = \sum_{i=0}^{k} z^i P_i, \quad P_i \in \F^{m\times n},
\end{equation}
of grade $k$, i.e., $P_k$ is allowed to be zero, and with one of the most influential families of pencils, i.e., matrix polynomials of grade $1$, associated to $P(z)$: the vector space of pencils $\mathbb{DL}(P)$ introduced by Mackey, Mackey, Mehl, and Mehrmann in \cite{MacMMM06} in the square case $m=n$. The matrix polynomial $P(z)$ is said to be regular if $m=n$ and its determinant is not identically zero. Otherwise, $P(z)$ is said to be singular. Regular and singular matrix polynomials arise in many applications \cite{nlevp,GohLR09,tisseur-meerbergen,VD-ieee-1981} and one of the most reliable methods for computing their eigenvalues, eigenvectors %(in the regular case)
 and minimal indices and bases %(in the singular case)
  is via {\em linearization}. See Section \ref{sec-background} for precise definitions of these and other concepts appearing in this introduction.

A linearization of $P(z)$ in \eqref{eq.intropoly} is a matrix pencil $L(z)$ with the same finite eigenvalues and associated partial multiplicities as $P(z)$ and also with the same number of right and left minimal indices. The linearization $L(z)$ is strong if, in addition, it has an eigenvalue at infinity whenever $P(z)$ has an eigenvalue at infinity, with the same partial multiplicities.  Linearizations and strong linearizations are very important in numerical computations because, when $\F \subseteq \mathbb{C}$, there exist reliable algorithms for approximating eigenvalues and computing minimal indices of matrix pencils, either regular \cite{moler-stewart-QZ} or singular \cite{VD-kron-1979}. Applying these algorithms to a strong linearization allows us to get the corresponding magnitudes of $P(z)$. As a consequence many families of strong linearizations of matrix polynomials have been developed and studied in the last two decades. We emphasize in this context the seminal work \cite{MacMMM06}. See also, for instance, \cite{antoniou-vol,DLPV17,hmmtsymlin,MacMMM06-2,NNT17} and the references therein among many other references on linearizations.

Among the many existing families of linearizations, the one based on the vector space of pencils $\mathbb{DL}(P)$ introduced in \cite{MacMMM06} has been particularly relevant and influential for several reasons in the study and numerical treatment of {\em regular} matrix polynomials $P(z)$. To begin with, for many classes of structured regular matrix polynomials $P(z)$ appearing in applications (symmetric, Hermitian, palindromic, odd, even, ...), the pencils in $\mathbb{DL}(P)$ allow us to construct linearizations with the same structure as $P(z)$ \cite{hmmtsymlin,MacMMM06-2}, which has important theoretical and numerical  consequences. In addition, the linearizations in $\mathbb{DL}(P)$ allow us to recover the left eigenvectors of $P(z)$ from the corresponding ones of the linearization exactly in the same way as the right eigenvectors \cite{hmmtsymlin,MacMMM06-2}. Moreover, the numerical properties of the linearizations in $\mathbb{DL}(P)$ were studied in \cite{higham-li-tisseur,higham-mack-tiss-cond}, where it was argued that, under some reasonable assumptions, some linearizations in $\mathbb{DL}(P)$ have satisfactory properties with respect to the conditioning of eigenvalues and the backward errors of eigenpairs, though some recent results \cite{bueno-perez} indicate that such conclusions should be reconsidered. Finally, we mention that the linearizations in $\mathbb{DL}(P)$ have played a fundamental role in the description and  parametrization of structure-preserving transformations for matrix polynomials in \cite{garvey-tisseur}, and that a generalization of the space $\mathbb{DL}(P)$, that goes beyond ansatz polynomials of bounded degree by considering ansatz polynomials of any degree and even more generally ansatz functions, was studied in \cite{NNT17}.

An important remark is that, given a {\em regular} matrix polynomial $P(z)$ of grade $k$, not every pencil in the vector space $\mathbb{DL}(P)$ is a linearization for $P(z)$. In order to explain this fact, let us recall that each nonzero pencil in $\mathbb{DL}(P)$ is uniquely determined by a nonzero vector $\omega = \begin{bmatrix} \omega_1 & \omega_2 & \cdots & \omega_k
\end{bmatrix}^T \in \F^k$, called \emph{ansatz vector} in \cite{MacMMM06}, which gives rise to a scalar polynomial $v(z) = \omega_1\,  z^{k-1} + \omega_2 \, z^{k-2} + \cdots + \omega_k$ of grade $k-1$, called \emph{ansatz polynomial} in \cite{NNT17}. It was proved in \cite[Theorem 6.7]{MacMMM06} that, for $P(z)$ regular, the pencil in $\mathbb{DL}(P)$ with ansatz polynomial $v(z)$ is a linearization for $P(z)$ if and only if the sets of eigenvalues of $P(z)$ and of roots of $v(z)$, each including possibly infinite eigenvalues or roots\footnote{Recall that the scalar polynomial $v(z) = \omega_1\,  z^{k-1} + \omega_2 \, z^{k-2} + \cdots + \omega_k$ of grade $k-1$ has a root at infinity of multiplicity $\ell < k-1$ if $\omega_1 = \cdots = \omega_\ell = 0$ and $\omega_{\ell + 1} \ne 0$.}, are disjoint. This neat result was termed in \cite{MacMMM06} the \emph{eigenvalue exclusion theorem} and led in \cite[Theorem 6.8]{MacMMM06} to the conclusion that almost all pencils in $\mathbb{DL}(P)$ are linearizations of $P(z)$ when $P(z)$ is regular. Moreover, it was proved in \cite[Theorem 4.3]{MacMMM06} that a pencil in $\mathbb{DL}(P)$ is a linearization for $P(z)$ if and only if it is a strong linearization for $P(z)$. On the other hand, neither fact remains true when the matrix polynomial $P(z)$ is singular, since it was proved in \cite[Theorem 6.1]{DDM09} that if $P(z)$ is a singular square matrix polynomial, then none of the pencils in $\mathbb{DL}(P)$ is a linearization for $P(z)$.

Theorem 6.1 in \cite{DDM09} is, at a first glance, a very discouraging result, since it seems to indicate that the good properties mentioned above of the pencils in $\mathbb{DL}(P)$ cannot be used at all for singular matrix polynomials. Moreover, it seems to put $\mathbb{DL}(P)$ at disadvantage with respect to other families of linearizations of matrix polynomials, which are linearizations for both regular and singular matrix polynomials \cite{bueno-curlett,das-bora,DDM09,DDM10,DDM12,DLPV17,nofer-perez}. The main purpose of this paper is to prove that the pencils in $\mathbb{DL}(P)$, though they are not linearizations, can still be used to compute the complete eigenstructure of singular matrix polynomials, since they reveal all the structural data of any singular matrix polynomial under exactly the same hypotheses of the eigenvalue exclusion theorem proved in \cite{MacMMM06} for regular polynomials. Thus, for instance, the pencils in $\mathbb{DL}(P)$ can be used to construct Hermitian pencils that reveal all the structural data of any, possibly singular, Hermitian matrix polynomial. In this context, the main result of this paper is the following theorem.
\begin{theorem}\label{thm:main}
Let $P(z)$ be an $m\times n$, possibly singular, matrix polynomial with coefficients in $\F^{m \times n}$ and of grade $k\geq 2$, right nullity $p$ and left nullity $q$, and let $v(z)$ be a scalar polynomial of grade $k-1$ with coefficients in $\F$. Suppose that the sets of eigenvalues of $P(z)$ and of roots of $v(z)$, each including possibly infinite eigenvalues or roots, are disjoint. Denote by $L(z)$ the pencil in $\mathbb{DL}(P)$ with ansatz polynomial $v(z)$. Then:
\begin{enumerate}
\item The right minimal indices of $P(z)$ are $\gamma_1 \leq \dots \leq \gamma_p$ if and only if the right minimal indices of $L(z)$ are $\underbrace{0= \dots= 0}_{p(k-1) \ \mathrm{times}} \leq \gamma_1 \leq \dots \leq \gamma_p$;
\item The left minimal indices of $P(z)$ are $\eta_1 \leq \dots \leq \eta_q$ if and only if the left minimal indices of $L(z)$ are $\underbrace{0 = \dots = 0}_{q(k-1) \ \mathrm{times}} \leq \eta_1 \leq \dots \leq \eta_q$;
\item $\lambda \in \F \cup \{ \infty \}$ is an eigenvalue of $P(z)$ if and only if $\lambda$ is an eigenvalue of $L(z)$, and its partial multiplicities as an eigenvalue of $P(z)$ and as an eigenvalue of $L(z)$ coincide.
\end{enumerate}
\end{theorem}
Note that, if $P(z)$ is regular, then Theorem \ref{thm:main} implies the eigenvalue exclusion theorem in \cite{MacMMM06}: indeed, items 1 and 2 become vacuous as there are no minimal indices (that is, in item 1 $p=0$ and in item 2 $q =0$), while item 3 is equivalent to saying that $L(z)$ is a strong linearization of $P(z)$. Hence, we can view Theorem \ref{thm:main} as a generalization of \cite[Theorem 6.7]{MacMMM06}. Observe that Theorem \ref{thm:main} immediately implies that, if $P(z)$ is singular, then $L(z)$ is not a linearization for $P(z)$, because the number of left (or right) minimal indices of $P(z)$ is different than the number of left (or right) minimal indices of $L(z)$ \cite{DDM14}. However, the key point is that for practical purposes such a pencil is not worse than a strong linearization: the minimal indices of $L(z)$ are those of $P(z)$ with the addition of a fixed number of extra zero minimal indices. Hence, the additional minimal indices can be trivially identified since the right and left nullities of $P(z)$ are just those of $L(z)$ divided by $k$. In fact, arguably the recovery of the minimal indices from the pencils in $\mathbb{DL}(P)$ for a singular $P(z)$ is simpler than their recovery from other families of pencils that are strong linearizations, including the classical companion pencils, whose minimal indices are not exactly those of $P(z)$, but shifted sequences of them \cite{bueno-curlett,DDM09,DDM10,DLPV17,nofer-perez}.

We disclose in advance that the proof of Theorem \ref{thm:main} requires considerable work and that it relies heavily on the bivariate polynomial approach to $\mathbb{DL}(P)$ pencils introduced in \cite{NNT17}, that maps the pencils in $\mathbb{DL}(P)$ into Lerer-Tismenetsky B\'ezoutians. We emphasize that this approach also allows us to consider for the first time in the literature $\mathbb{DL}(P)$ pencils for matrix polynomials $P(z)$ that may be rectangular without extra effort.

The results in this paper are connected with those in \cite{DQVD-selfconjugate} concerning matrix polynomials, though the approaches in both papers are completely different. Given an $m\times n$ complex matrix polynomial, the authors of \cite{DQVD-selfconjugate} consider only one pencil $L(z)$ in $\mathbb{DL} (P)$, the one with ansatz polynomial $v(z) = 1$, i.e., the ansatz polynomial with $\infty$ as unique root. For this particular pencil $L(z)$, it was proved in \cite{DQVD-selfconjugate}, without assuming the eigenvalue exclusion condition, that Theorem \ref{thm:main} holds with the exception of item 3 for the eigenvalue $\lambda = \infty$, for which it was proved that though it has different partial multiplicities in $L(z)$ and in $P(z)$, the ones in $P(z)$ can be recovered from the ones in $L(z)$. The approach in \cite{DQVD-selfconjugate} is \emph{ad hoc} in that sense that it relies on making some matrix manipulations that exploit the very peculiar structure of the pencil in $\mathbb{DL} (P)$ with ansatz polynomial $v(z) = 1$ and on the concept of strongly minimal linearization of a matrix polynomial. The latter is different from the standard definitions of linearization and strong linearization coming from \cite{GohLR09}. Hence, the results in \cite{DQVD-selfconjugate} are presented in a form rather different from Theorem \ref{thm:main}. In this paper, in contrast, we systematically consider all the (uncountably many) pencils in the vector space $\mathbb{DL} (P)$  and we prove that they satisfy Theorem \ref{thm:main} working with their representation as B\'ezoutians via a bivariate polynomial mapping.

It is worth mentioning that, with the exception of the particular pencil considered in \cite{DQVD-selfconjugate}, this paper presents the first example of a wide family of pencils related to a matrix polynomial $P(z)$ that are not linearizations, according to the standard definition introduced by Gohberg, Lancaster and Rodman \cite{GohLR09}, but that still allows for the recovery of all the structural data of $P(z)$. Since the recovery of such structural data is what is really wanted in applications, the results we have obtained for the pencils in $\mathbb{DL} (P)$ indicate that the Gohberg, Lancaster and Rodman's concepts of linearization and strong linearization may be unnecessarily rigid.

The paper is organized as follows. Section \ref{sec-background} summarizes the background material that is necessary in the rest of the paper. Some preliminary technical results that are used in the proofs of the main results are presented and proved in Section \ref{sec-preliminary}. Items 1 and 2 of Theorem \ref{thm:main} are proved in Section \ref{sec-minimalbasisdlp}, while item 3 is proved in Section \ref{sec-elemdivisors}.  Section \ref{sec.recovery} studies how to recover minimal bases, eigenvectors, and root polynomials of a possibly singular matrix polynomial $P(z)$ from those of pencils in the vector space $\mathbb{DL}(P)$. Finally, some conclusions are presented in Section \ref{sec-conclusions}.

 \section{Background material} \label{sec-background}
 In this section, we revise the necessary background material. Throughout this paper, we consider an algebraically closed field $\F$ with characteristic $0$, the ring of polynomials $\F [z]$ with coefficients in $\F$, and its field of fractions $\F (z)$. The sets of $m \times n$ constant matrices over $\F$, of $m \times n$ matrix polynomials over $\F$, and of $m \times n$ rational matrices over $\F$ are denoted by $\F^{m \times n}$, $\F[z]^{m \times n}$ and $\F(z)^{m \times n}$, respectively.

 \subsection{Basics on matrix polynomials}
 The normal rank, or simply the rank, of a matrix polynomial $P(z) \in \F[z]^{m \times n}$ is denoted by $\mbox{rank} \, P(z)$ and is the size of its largest non-identically-zero minor or, equivalently, its rank over the field $F(z)$. The matrix polynomial $P(z)$ is regular if it is square and $\det P(z)$ is not identically zero. Otherwise, $P(z)$ is said to be singular. An element $\lambda \in \F$ is a finite eigenvalue of $P(z)$ if the rank of the constant matrix $P(\lambda) \in \F^{m\times n}$ is less than the normal rank of $P(z)$, i.e., if $\mbox{rank} \, P(\lambda) < \mbox{rank} \, P(z)$.

 A square matrix polynomial $U(z) \in \F[z]^{n \times n}$ is said to be unimodular if its determinant is a nonzero element of $\F$ or, equivalently, if it is invertible and its inverse is a matrix polynomial.

The Smith canonical form of $P(z)$ with $r = \mbox{rank} \, P(z)$ is $S(z)=U(z) P(z) V(z)$, where $U(z),V(z)$ are unimodular matrix polynomials and $S(z) = \diag (d_1 (z), d_2 (z), \ldots , d_r (z), 0, \ldots, 0)$ is diagonal and such that each nonzero diagonal entry $d_1 (z), d_2 (z), \ldots , d_r (z)$ divides the next one \cite{gantmacher}. The diagonal elements of $S(z)$ are unique up to multiplication by units of $\F[z]$ and the nonzero ones are called invariant factors of $P(z)$. Given $\lambda \in \F$, each  invariant factor can be written uniquely as $d_i (\lambda) = (z-\lambda)^{\ell_i} \, \omega_i(z)$, where $\omega_i(\lambda)\neq 0$ and $\ell_i$ is a nonnegative integer, for $i= 1, \ldots , r$. Thus, $\lambda \in \F$ is an eigenvalue of $P(z)$ if and only if at least one $\ell_i$ is strictly larger than zero. Those $\ell_i$ that are positive are called the partial multiplicities of $\lambda$ as an eigenvalue of $P(z)$ and the corresponding factors $(z-\lambda)^{\ell_i}$ are the elementary divisors of $P(z)$ for the eigenvalue $\lambda$.

The partial multiplicities at infinity of $P(z)$ are defined after fixing a grade $k$ for $P(z)$, where $k$ is an integer greater than or equal to the degree of $P(z)$, denoted by $\deg P(z)$ and defined as the largest degree of its entries. Then, the reversal of $P(z)$ with respect to the grade $k$ is defined as $\Rev_k P(z) = z^k P(z^{-1})$. Very often the grade $k$ is chosen to be equal to $\deg P(z)$; if this is the case we will omit the suffix $k$ and simply write $\Rev P(z)$. The partial multiplicities of $\infty$ in $P(z)$, seen as a grade $k$ polynomial, are then defined to be equal to the partial multiplicities of $0$ for $\Rev_k P(z)$. It is clear that different choices of the grade $k$ lead to different partial multiplicities at $\infty$: more precisely, if $k = \deg P(z)+h$, for some $h>0$, and $\rank P(z)=r$, then $\Rev P(z)$ has partial multiplicities at $0$ equal to $0 < \eta_1 \leq \dots \leq \eta_s$ if and only if $\Rev_k P(z)$ has partial multiplicities at $0$ equal to $\underbrace{h=\dots=h}_{r-s \ \mathrm{times}} \leq \eta_1 +h \leq \dots \leq \eta_s+h$. Note that here $s$ may be or may be not zero, and accordingly $\infty$ may be not or may be an eigenvalue of $P(z)$ when $k = \deg P(z)$. However, if $k > \deg P(z)$ and $r>0$, then $\infty$ is always an eigenvalue of $P(z)$ with $r$ partial multiplicities.

The right kernel, or right null space, of $P(z) \in \F[z]^{m \times n}$ over the field $\F (z)$ is defined as
$$
 \ker P(z) = \{ x(z) \in \F (z)^{n \times 1} \, : \, P(z) x(z) = 0\}\, .
$$
The left kernel of $P(z)$ over $\F (z)$ is defined as the right kernel of the transpose of $P(z)$, which is denoted as $P(z)^T$. In this paper, we mainly work with right kernels, since the results for left kernels are completely analogous, and for this reason we will for simplicity just write ``kernel'' to mean ``right kernel''.
It is clear that $\ker P(z)$ has polynomial bases, i.e., bases whose vectors belong to $\F [z]^{n \times 1}$. Among the polynomial bases of $\ker P(z)$, those which have minimal sum of the degrees of their vectors are called (right) minimal bases of $P(z)$. There are infinitely many (right) minimal bases of $P(z)$, but all of them have the same ordered sequence of the degrees of their vectors \cite{For75}. Such degrees are called the right minimal indices of $P(z)$. Left minimal bases and indices of $P(z)$ are defined as the right ones of $P(z)^T$.

The index sum theorem \cite{DDM14} states that, for every matrix polynomial $P(z) \in \F[z]^{m \times n}$ of grade $k$ and rank $r$, the sum of all the left and right minimal indices plus the sum of all finite and infinite partial multiplicities is equal to $k\, r$.

\subsection{Linearizations of matrix polynomials and $\mathbb{DL}(P)$ pencils}
In applications of matrix polynomials, one is usually interested in computing their eigenvalues, together with their associated partial multiplicities, and their minimal indices. Often, this is done by linearizing the matrix polynomial $P(z) \in \F [z]^{m \times n}$ one is interested in.  Linearizations are pencils that allow us to obtain the required information of $P(z)$ by coupling a theoretical analysis that relates this information with its analogue for $L(z)$ together with an algorithm such as, for instance, those in \cite{moler-stewart-QZ,VD-kron-1979}. More precisely a pencil $L(z)=L_0 + L_1 z$ is a linearization of $P(z)$ if there exist unimodular matrix polynomials $U(z),V(z)$ such that $U(z)L(z)V(z)=I \oplus P(z)$; if, moreover, $\Rev_1 L(z)$ is a linearization for $\Rev_k P(z)$, then $L(z)$ is said to be a strong linearization for $P(z)$ considered as a polynomial of grade $k$ \cite{GohLR09}. Strong linearizations preserve all the finite and infinite eigenvalues, together with their partial multiplicities, as well as the dimensions of the left and right kernels. However, they typically do not preserve minimal indices \cite{DDM14}, though the minimal indices of the most standard classes of strong linearizations available in the literature allow us to obtain very easily those of the matrix polynomial \cite{DDM09,DLPV17}.

In this paper, we will study a family of pencils associated with a matrix polynomial $P(z)\in \F [z]^{m \times n}$ that are not linearizations according to the definition above (in particular, they do not preserve the dimensions of the left and right kernels), but that still preserve all the finite and infinite eigenvalues of $P(z)$, together with their partial multiplicities, and moreover allow for the recovery of the minimal indices of $P(z)$ in a way which is (at least) as simple as for the most standard types of linearizations. This family of pencils is described in the rest of this subsection. We emphasize that in this description the considered matrix polynomial is arbitrary, i.e., it may be rectangular, square and singular, or square and regular. This generality is in contrast with previous references \cite{MacMMM06,NNT17} that deal with this family but restrict their attention to square matrix polynomials that, moreover, are often also required to be regular. It is easy to check that those properties of the pencils in this family that are needed in this work hold for general matrix polynomials, which motivates us to consider in this work fully general matrix polynomials.

The so-called $\mathbb{DL}(P)$ vector space of pencils associated with a matrix polynomial $P(z)\in \F [z]^{m \times n}$ was first defined by Mackey, Mackey, Mehl and Mehrmann in \cite{MacMMM06}. For a fixed matrix polynomial $P(z)$ of grade $k$ and size $m\times n$, each pencil in $\mathbb{DL}(P)$ has size $k m\times kn$ and is uniquely determined by a nonzero ansatz vector $\omega$ in $\F^k$. In turn, in this paper we identify nonzero ansatz vectors with nonzero scalar ansatz polynomials\footnote{In \cite{MacMMM06}, what we call the ansatz polynomial was instead called the $v$-polynomial. However, we prefer to follow instead the nomenclature of \cite{NNT17} because it fits better with the B\'{e}zoutian description of $\mathbb{DL}(P,v)$, and the latter is crucial to present relatively simple proofs of our results.} of grade $k-1$ via the bijection
 \begin{equation} \label{eq.omegatov}   \omega \in \F^k \mapsto v(x)=\omega^T \begin{bmatrix}
x^{k-1}\\
\vdots\\
x\\
1
 \end{bmatrix} \in \F[x]_{k-1},\end{equation}
where $\F[x]_{k-1}$ denotes the vector space of scalar polynomials in the variable $x$ of grade $k-1$, or, equivalently, of degree at most $k-1$.
If $L(z)$ is the pencil in $\mathbb{DL}(P)$ associated with the ansatz polynomial $v(z)$, we denote it by $L(z)=:\mathbb{DL}(P,v)$, omitting in the right-hand side of this identity the variable $z$ for brevity. The eigenvalue exclusion theorem \cite[Theorem 6.7]{MacMMM06} states that, when $P(z)$ is a regular matrix polynomial, $\mathbb{DL}(P,v)$ is a strong linearization for $P(z)$ if and only if the set of the roots of $v(z)$ (when seen as a grade $k-1$ scalar polynomial) and the set of the eigenvalues of $P(z)$ (when seen as a grade $k$ matrix polynomial) are disjoint. For a given regular matrix polynomial $P(z)$, this ``disjointness'' condition is satisfied for almost all ansatz polynomials $v(z)$, which motivates to label $\mathbb{DL}(P)$ as a vector space of ``potential linearizations'' for $P(z)$. In contrast, it was proved in \cite[Theorem 6.1]{DDM09} that if $P(z)$ is square and singular, then none of the pencils in $\mathbb{DL}(P)$ is a linearization for $P(z)$.

In \cite{NNT17}, Nakatsukasa, Noferini and Townsend proposed a convenient interpretation of block matrices as bivariate matrix polynomials. Any $k \times k$ block matrix $B$ of size $km \times kn$ with blocks $B_{ij}$ each of size $m \times n$, $1\leq i,j \leq k$, is associated with the following $m\times n$ bivariate matrix polynomial of grade $k - 1$ in both variables
 \begin{equation}\label{eq:bivariatebijection}
  F(x,y) = \begin{bmatrix}
 y^{k-1} I_m & \cdots & y I_m & I_m
 \end{bmatrix} B \begin{bmatrix}
 x^{k-1} I_n\\
 \vdots\\
 x I_n\\
 I_n
 \end{bmatrix}  = \sum_{i=0}^{k-1} \sum_{j=0}^{k-1} y^i x^j B_{k-i,k-j} ,
\end{equation}
where $I_m$ and $I_n$ denote the identity matrices of sizes $m\times m$ and $n\times n$, respectively.
 This map is manifestly a bijection between $m \times n$ bivariate matrix polynomials of grade $k-1$ in both $x$ and $y$ and $km \times kn$ matrices partitioned into $k\times k$ blocks all of size $m\times n$. This framework has the advantage of greatly simplifying many operations that are associated with the pencils in the $\mathbb{DL}(P)$ vector space: see, for example, \cite[Tables 1 and 2]{NNT17}. In particular, it was shown in \cite{NNT17} that when $m=n$ the pencil $\mathbb{DL}(P,v)$ is mapped by the bijection above to the Lerer-Tismenetsky B\'{e}zoutian \cite{AJ76,LT82,LT86} associated with $P(x)$ and $(x-z)v(x)I_n$. This statement will be very important in this paper and deserves a detailed explanation in terms of bivariate polynomials and block matrices; moreover, we will generalize it to the case $m \neq n$. Note that given $P(z)\in \F [z]^{m \times n}$ of grade $k$, the pencil $\mathbb{DL}(P,v)=:L_1 z + L_0 \in \F [z]^{km \times kn}$ can be seen as a $k\times k$ block pencil or, equivalently, as a block matrix whose blocks belong to $\F[z]^{m \times n}$ and have grade $1$. Hence, we can view it via the bijection \eqref{eq:bivariatebijection} as a trivariate matrix polynomial in $x,y,z$ having grade $1$ in $z$ and grade $k-1$ in both $x$ and $y$. We can express concisely this Lerer-Tismenetsky-B\'{e}zout connection as follows: functionally, the bijection \eqref{eq:bivariatebijection} maps
\begin{equation}\label{eq:bezout}
 \mathbb{DL}(P,v) \in \F[z]^{km \times kn} \mapsto  \mathcal{B}_{P,v} (x,y,z) := \frac{P(y)(x-z)v(x)-P(x)(y-z)v(y)}{x-y} \in \F[x,y,z]^{m \times n},
\end{equation}
where we observe that $\mathcal{B}_{P,v} (x,y,z) = \mathcal{B}_{P,v} (y,x,z)$.
The relation \eqref{eq:bezout} was crucial in \cite{NNT17} to enormously simplify the original proof of the eigenvalue exclusion theorem when $P(z)$ is regular, and to derive some novel results. As said in the introduction (see Theorem \ref{thm:main}), in this paper we essentially aim to extend the eigenvalue exclusion theorem to the case where $P(z)$ is singular (square or rectangular). This extension means that if the set of the roots of $v(z)$ (when seen as a grade $k-1$ scalar polynomial) and the set of the eigenvalues of $P(z)$ (when seen as a grade $k$ matrix polynomial) are disjoint, then we will prove that $\mathbb{DL}(P,v)$ has the same eigenvalues and partial multiplicities as $P(z)$, while the minimal indices of $\mathbb{DL}(P,v)$ are precisely those of $P(z)$ together with some extra minimal indices equal to zero. Therefore, it is not surprising that the bijection \eqref{eq:bezout} plays also a central role in this paper. In this context, we emphasize that \eqref{eq:bezout} was derived in \cite{NNT17} without any need to assume that $P(z)$ is regular \cite[Remark 3.2]{NNT17}, though in many of the results that followed next in \cite{NNT17} it was assumed that $P(z)$ is regular.

\begin{remark}
 In principle, the bijection \eqref{eq:bivariatebijection} could be defined using \emph{any} polynomial basis for the space of scalar polynomials of degree at most $k-1$, i.e., the restriction to the monomial basis is convenient, but not necessary. This was indeed already observed in \cite{NNT17} and rediscovered in \cite{FS17}. It leads to an immediate generalization of the space $\mathbb{DL}(P)$. While we restrict to the monomial basis for the sake of a much simpler exposition, we observe that the results of this paper also generalize to $\mathbb{DL}(P)$ in non-monomial bases. This follows easily by observing that (a) a bijection between such a space and the ``traditional'' $\mathbb{DL}(P)$ in the monomial basis is given by multiplying each pencil in the generalized $\mathbb{DL}(P)$ by $M^T \otimes I_m$ on the left and by $M \otimes I_n$ on the right, where $M$ is the $k\times k$ nonsingular constant change-of-polynomial-basis matrix, and that (b) all the relevant structures, i.e., eigenvalues, partial multiplicities, minimal indices and bases, root polynomials, are well behaved under strict equivalence. This remark could have practical relevance, because $\mathbb{DL}(P)$ in certain non-monomial bases may be useful to exploit the better numerical properties of some basis. One example, albeit limited to the scalar case $m=n=1$, is presented in \cite{NNT15} where a method based on the B\'{e}zout resultant matrix (i.e., $\mathbb{DL}(P)$) in the Chebyshev basis was used to compute real common zeros of bivariate functions.
\end{remark}

\subsection{Root polynomials} \label{subsec.rootpolys}
Maximal sets of root polynomials \cite{DopicoNoferini, Nof12} (see also \cite{GohLR09} for the regular case and \cite{NVD23} for an extension to rational functions) are another crucial tool in this paper. They are useful for the analysis of singular matrix polynomials because they carry the relevant information about the partial multiplicities of a matrix polynomial at a given eigenvalue $\lambda$. More precisely, suppose that $M(z) \in \F[z]^{n \times p}$ is a right minimal basis for the matrix polynomial $P(z) \in \F[z]^{m \times n}$, which is a short way of stating that the columns of $M(z)$ form a minimal basis of the right kernel of $P(z)$. Then, for any $\lambda \in \F$ we define $\ker_\lambda P(z)$ as the subspace of $\F^n$ spanned by the columns of $M(\lambda)$. A polynomial vector $r(z) \in \F[z]^n$ is said to be a \emph{root polynomial} at $\lambda$ of  order $\ell \geq 1$ for $P(z)$ if
\begin{enumerate}
\item $P(z)\, r(z) = (z-\lambda)^\ell w(z)$, for some $w(z) \in \F[z]^m$ such that $w(\lambda)\neq 0$, and
\item $r(\lambda) \not \in \ker_\lambda P(z)$.
\end{enumerate}
Suppose that $\{ r_i(z) \}_{i=1}^t$ is a set of root polynomials at $\lambda$ for $P(z)$, with orders $\ell_1 \geq \dots \geq \ell_t$, and let $R(z) \in \F[z]^{n \times t}$ be the matrix whose $i$th column is $r_i(z)$. The set $\{ r_i(z) \}_{i=1}^t$  is said to be $\lambda$-independent if the matrix $\begin{bmatrix}
M(\lambda) & R(\lambda)
\end{bmatrix}$ has full column rank; it is said to be complete if it is $\lambda$-independent and $t = \dim \ker P(\lambda) - \dim \ker_\lambda P(z)$; and it is said to be maximal if it is complete and satisfies the property that, for each $i=1,\dots,t$, there is no root polynomial $\hat r(z)$ such that $r_1(z),\dots,r_{i-1}(z),\hat r(z)$ are $\lambda$-independent and the order of $\hat r(z)$ is strictly larger than the order of $r_i(z)$. It was proved in \cite{DopicoNoferini,NVD23} that all the concepts defined above are independent of the particular right minimal basis $M(z)$ of $P(z)$ that is considered and that root polynomials at $\lambda$ for $P(z)$ exist if and only if $\lambda$ is a finite eigenvalue of $P(z)$.

All maximal sets of root polynomials of a given matrix polynomial $P(z)$ at a given point $\lambda$ have the same orders \cite[Theorem 4.1-Item 2]{DopicoNoferini} and these orders coincide with the partial multiplicities of $\lambda$ as an eigenvalue of $P(z)$ \cite[Theorem 4.2]{DopicoNoferini}. Moreover, the following theorem that appeared in \cite[Theorem 4.1-Item 3]{DopicoNoferini} will be important in this paper.

\begin{theorem}\label{whenismaximal}
Let $\{ r_i(z) \}_{i=1}^t$ be a complete set of root polynomials at $\lambda \in \F$ for $P(z) \in \F[z]^{m \times n}$ with orders $\kappa_1 \geq \dots \geq \kappa_t$. Suppose, moreover, that the partial multiplicities of $\lambda$ as an eigenvalue of $P(z)$ are $\ell_1 \geq \dots \geq \ell_t$. Then, $\ell_i \geq \kappa_i$ for all $i=1,\dots,t$. Moreover, the following are equivalent:
\begin{enumerate}
\item $\{ r_i(z) \}_{i=1}^t$ is a maximal set of root polynomials at $\lambda$ for $P(z)$;
\item $\ell_i = \kappa_i$ for all $i=1,\dots,t$;
\item $\sum_{i=1}^t \kappa_i = \sum_{i=1}^t \ell_i$.
\end{enumerate}
\end{theorem}

The definitions above of root polynomials and $\lambda$-independent, complete and maximal sets of root polynomials correspond to the right versions of these concepts. Left root polynomials and $\lambda$-independent, complete and maximal sets of left root polynomials are defined as the right ones for $P(z)^T$. The definition of root polynomials at infinity for $P(z)$ requires certain care, but their properties are analogous to those of root polynomials at finite points. They can be found in \cite[Section 6]{DopicoNoferini}.

\section{Preliminary results} \label{sec-preliminary} This section includes some auxiliary results that are useful to prove the main results in this paper. The lemmata in  Subsection \ref{subsec-vander} are essential to study in Section \ref{sec-minimalbasisdlp} the minimal indices of the pencils in $\mathbb{DL} (P)$. The results in Subsection \ref{subsec-mobiusdlp} allow us to simplify many of the proofs in the rest of the paper, because they show, via a M\"{o}bius transformation, that it is enough to establish them only for matrix polynomials $P(z)$ and ansatz polynomials $v(z)$ in \eqref{eq:bezout} without eigenvalues and roots at infinity, respectively.

\subsection{Technical lemmata on certain Vandermonde-structured matrices} \label{subsec-vander}
In this subsection we provide some results on matrices that have a tensor product Vandermonde structure; these results are crucial for the construction in Section \ref{sec-minimalbasisdlp} of a minimal basis of $\ker L(z)$, where $L(z)=\mathbb{DL}(P,v)$, for a singular matrix polynomial $P(z)$ and an ansatz polynomial $v(z)$ whose roots are not eigenvalues of $P(z$).

Let us start with Lemma \ref{important}, that assumes a tensor product structure with a classical (i.e. not confluent) Vandermonde matrix.

\begin{lemma}\label{important} Let $\ell \leq k$ be two positive integers. For $i=1,\dots,\ell$, let $A_i \in \mathbb{F}^{n \times m_i}$ with $\rank A_i = m_i$  and let $\alpha_i \in \mathbb{F}$ be such that $\alpha_i \neq \alpha_j$ if $i \neq j$. Let
$$V(z)=\begin{bmatrix}
z^{k-1}\\
\vdots \\
z\\
1
\end{bmatrix}$$ be the Vandermonde vector of degree $k-1$ depending on the variable $z \in \F$. Then the matrix $C=\begin{bmatrix}
V(\alpha_1) \otimes A_1 & \cdots & V(\alpha_\ell) \otimes A_{\ell}
\end{bmatrix} \in \mathbb{F}^{k n \times m}$ has full column rank $m=\sum_{i=1}^{\ell} m_i$.
\end{lemma}

\begin{proof}
Let $S=\begin{bmatrix}
V(\alpha_1) & \cdots & V(\alpha_\ell)
\end{bmatrix}$, which can be seen as a $k \times \ell$ submatrix of a $k\times k$ Vandermonde matrix. Then, standard properties of Vandermonde matrices imply that $S$ has full column rank $\ell$. Since we can factorize $C=(S \otimes I_n) \diag(A_1,\dots,A_\ell)$, we have $\rank C=\rank \diag(A_1,\dots,A_\ell)=m$.
\end{proof}

We have stated Lemma \ref{important} mainly for illustrative purposes, as it is simpler to analyze than the general case, but to our goals in Section \ref{sec-minimalbasisdlp} it is only useful if all the roots of the ansatz polynomial $v(z)$ are simple. This is in fact not typical in practice since the most common choices for $v(z)$ are either $1$ or $z^{k-1}$ \cite{DQVD-selfconjugate,hmmtsymlin,MacMMM06}, that have a unique root of multiplicity $k-1$ at $\infty$ and at $0$, respectively. However, by looking at the proof, we observe that Lemma \ref{important} would still hold if we replaced $S$ by any full column rank matrix, not necessarily Vandermonde. Our next goal is to generalize in Lemma \ref{confvanderparttwo} the result to a tensor product structure inherited by a \emph{confluent} Vandermonde matrix (see \cite[Chapter 22]{higham-error} for information on confluent Vandermonde matrices). Lemma \ref{confvanderpartone} is a preliminary step for this goal.

\begin{lemma}\label{confvanderpartone}
Let $\ell \leq k$ be two positive integers, let $\alpha \in \F$, let $M(z) \in \F[z]^{n \times p}$ be a matrix polynomial such that $M(\alpha)$ has full column rank, and let $V(z) \in \F[z]^{k \times 1}$ be defined as in Lemma \ref{important}. Let us denote by $[V(z) \otimes M(z)]'$ and $[V(z) \otimes M(z)]^{(\beta)}$ the derivative and the $\beta$th derivative with respect to $z$, respectively, and by $[\, \cdot \, ]|_{z=\alpha}$ evaluation at $\alpha$. Then the matrix
\[  M_{\alpha,k,\ell,M(z)} =  \left. \begin{bmatrix}
V(z) \otimes M(z) & [V(z) \otimes M(z)]' & \cdots & [V(z) \otimes M(z)]^{(\ell-1)}
\end{bmatrix} \right|_{z=\alpha} \in \F^{k n \times \ell p}  \]
has full column rank.
\end{lemma}

\begin{proof}
The block entry of $M_{\alpha,k,\ell,M(z)}$ (before evaluation at $\alpha$) in block row $k-\phi$ and block column $\beta+1$ is
\[   [z^\phi M(z)]^{(\beta)} = \sum_{j=0}^\beta \, [z^\phi]^{(j)}  \begin{pmatrix}
\beta \\
j
\end{pmatrix} [M(z)]^{(\beta-j)}.    \]
Since the $j$th derivative of $z^\phi$ is $0$ for $j > \phi$ and in the summation only the first $\beta$ derivatives of $M(z)$ appear, we can encode these relations by expressing $M_{\alpha,k,\ell,M(z)}$ as the product $M_{\alpha,k,\ell,M(z)}=AB$, where
\begin{equation}\label{eq:structureofA}
A = \begin{bmatrix}
V^{(\ell-1)}(\alpha) &  \dots  & V'(\alpha) & V(\alpha)
\end{bmatrix}  \otimes I_n   \in \F^{ kn \times \ell n}
\end{equation}
corresponds to the last columns of a block upper triangular matrix, and
\begin{equation}\label{eq:structureofB}
B = \begin{bmatrix}
0 & 0 &  & \dots & M(\alpha)\\
\vdots & \vdots & & & \vdots\\
0 & 0 & M(\alpha) & \dots & \begin{pmatrix}
\ell-1\\
2
\end{pmatrix} M^{(\ell-3)}(\alpha)\\
0 & M(\alpha) & 2 M'(\alpha) & \dots & (\ell-1) M^{(\ell-2)}(\alpha)\\
M(\alpha) & M'(\alpha) & M''(\alpha) & \dots & M^{(\ell-1)}(\alpha)
\end{bmatrix} \in \F^{\ell n \times \ell p},
\end{equation}
is block lower antitriangular. Observe that, generally, the block $B_{ij}$ is $\begin{pmatrix}
j-1\\
\ell-i
\end{pmatrix} M^{(i+j-\ell-1)}(\alpha)$ if $i+j > \ell$ and $0$ otherwise.
Note that $A$ is a $k\times \ell$ submatrix of a $k\times k$ confluent Vandermonde matrix tensor product the identity, and hence it has full column rank. As a consequence, $\rank M_{\alpha,k,\ell,M(z)} = \rank B = \ell \rank M(\alpha) = \ell p$.
\end{proof}

Lemma \ref{confvanderpartone} implies the announced generalization of Lemma \ref{important}.

\begin{lemma}\label{confvanderparttwo}
Let $\ell_1, \ldots , \ell_s , k$ be positive integers such that $\ell:=\sum_{i=1}^s \ell_i \leq k$ and let $V(z) \in \F[z]^{k \times 1}$ be defined as in Lemma \ref{important}. For $i=1,\dots, s$, let $\alpha_i \in \mathbb{F}$  be such that $\alpha_i \neq \alpha_j$ if $i \neq j$, let $M_i(z) \in \F[z]^{n \times p_i}$ be a matrix polynomial such that $M_i(\alpha_i)$ has full column rank, and let $M_{\alpha_i,k,\ell_i,M_i(z)} \in \F^{k n \times \ell_i p_i} $ be defined as in Lemma \ref{confvanderpartone}. Then
\[ C = \begin{bmatrix}
M_{\alpha_1,k,\ell_1,M_1(z)} & M_{\alpha_2,k,\ell_2,M_2(z)} & \cdots & M_{\alpha_s,k,\ell_s,M_s(z)}
\end{bmatrix} \in \, \F^{k n \times \left( \sum_{i=1}^{s} \, \ell_i p_i \right) }\]
has full column rank.
\end{lemma}

Before proving Lemma \ref{confvanderparttwo}, we observe that $\ell_i$ and $M_i(z)$ are allowed to differ for different values of $i$ and that the matrices $M_i(z)$ may also have different number of columns (but they must have the same number of rows). However, it is a necessary assumption that $M_i(\alpha_i)$ has full column rank for all $i$. For our scopes in this paper, however, $M_i(z)$ will not depend on the index $i$ (whereas $\ell_i$ may). Observe also that $\sum_{i=1}^{s} \, \ell_i p_i \leq (\sum_{i=1}^{s} \, \ell_i) \max_i p_i \leq k n$. Thus, the number of columns of $C$ is less than or equal to the number of its rows.

\begin{proof}[Proof of Lemma \ref{confvanderparttwo}]
For each $i$, let $M_{\alpha_i,k,\ell_i,M_i(z)} = A_i B_i$ be a factorization as in the proof of Lemma \ref{confvanderpartone}. Then, $C=AB$ with
\[ A = \begin{bmatrix}
A_1 & A_2 & \cdots & A_s
\end{bmatrix} \in \F^{k n \times \ell n}, \qquad B = B_1 \oplus B_2 \oplus \dots \oplus B_s \in \F^{\ell n \times \left( \sum_{i=1}^{s} \, \ell_i p_i \right)},\]
where $A_i$, $B_i$ have the structures of \eqref{eq:structureofA} and \eqref{eq:structureofB}, respectively. Then, $A$ is a $k \times \ell$ submatrix of a $k\times k$, possibly confluent, Vandermonde matrix tensor product the identity. Therefore, it has full colum rank. Hence, $\rank C = \rank B = \sum_{i=1}^s \rank B_i = \sum_{i=1}^{s} \, \ell_i p_i$, proving the statement.
\end{proof}

\subsection{A commuting diagram between $\mathbb{DL}(P)$ pencils and M\"{o}bius transformations} \label{subsec-mobiusdlp}

Let $\gcd$ stand for greatest common divisor. Recall that, given a matrix polynomial $P(z)\in \F[z]^{m \times n}$ of degree $k$ and a M\"{o}bius function $r(z)=n(z)/d(z)$, where $n(z)=\alpha z + \beta$, $d(z)=\gamma z + \delta$ satisfy $\gcd(n(z),d(z))=1$ (or equivalently $A = \left[ \begin{smallmatrix}
\alpha & \beta \\
\gamma & \delta
\end{smallmatrix}\right]$
is invertible), the associated M\"{o}bius transformation of grade $g \geq k$ of $P(z)$ is
\[ \mathcal{M}_{g,r(z)}(P) = d(z)^g P\left( \frac{n(z)}{d(z)}   \right) \in \F[z]^{m \times n},\]
where we omit the dependence on $z$ in the left-hand side for simplicity.
These transformations have been studied in papers such as \cite{MacMMM14,Nof12}. Crucially, the map $z \mapsto n(z)/d(z)$ is an automorphism of $\F \cup \{ \infty \}$ (upon agreeing that, conventionally, for all nonzero $\alpha \in \F$, $\alpha/0=\infty$; moreover $n(\infty)=\alpha$ and $d(\infty)=\gamma$). In particular, the following results are known (see, for example, \cite[Theorems 5.3 and 7.5]{MacMMM14} and \cite[Theorem 4.1]{Nof12}).
\begin{theorem}\label{thm:eigenstrucutreafterMobius}
Suppose that $P(z) \in \F[z]^{m \times n}$ has grade $g$, that $Q(z)=\mathcal{M}_{g,r(z)}(P) \in \F[z]^{m \times n}$ is the M\"{o}bius transform of grade $g$ induced by $r(z)=n(z)/d(z)$ of $P(z)$, and assume that $Q(z)$ has also grade $g$. Let $\lambda,\mu \in \F \cup \{ \infty \}$ satisfy $\lambda  d(\mu)= n(\mu)$. Then:
\begin{enumerate}
\item $\lambda$ is an eigenvalue of $P(z)$ with nonzero partial multiplicities $\delta_1 \leq \dots \leq \delta_s$ if and only if $\mu$ is an eigenvalue of $Q(z)$ with nonzero partial multiplicities $\delta_1 \leq \cdots \leq \delta_s$;
\item $m_1(z),\ldots,m_p(z) \in \F[z]^{n \times 1}$ are a right minimal basis of $P(z)$ with right minimal indices $\nu_1 \leq \cdots \leq \nu_p$ if and only if $u_1(z),\ldots,u_p(z) \in \F[z]^{n \times 1}$ are a right minimal basis of $Q(z)$ with right minimal indices $\nu_1 \leq \cdots \leq \nu_p$, where $u_i(z)=\mathcal{M}_{\nu_i , r(z)}(m_i)$ for all $i=1,\ldots,p$.
\item $n_1(z),\ldots,n_q(z) \in \F[z]^{m \times 1}$ are a left minimal basis of $P(z)$ with left minimal indices $\eta_1 \leq \dots \leq \eta_q$ if and only if $v_1(z),\ldots,v_q(z) \in \F[z]^{m \times 1}$ are a left minimal basis of $Q(z)$ with left minimal indices $\eta_1 \leq \dots \leq \eta_q$, where $v_i(z)=\mathcal{M}_{\eta_i, r(z)}(n_i)$ for all $i=1,\dots,q$.
\end{enumerate}
\end{theorem}

With this in mind, we show in Theorem \ref{thm:infinitoacasa} how the operations of constructing a $\mathbb{DL}(P)$ pencil of a matrix polynomial and applying M\"{o}bius transformations relate. We prove the simple Lemma \ref{lemm:changebasemob} as a preliminary step.
\begin{lemma} \label{lemm:changebasemob} Let $n(z)=\alpha z + \beta , d(z)=\gamma z + \delta \in \F[z]$ be scalar polynomials of grade $1$ such that $\gcd(n(z),d(z))=1$ and let $\F[z]_{k-1}$ be the vector space of scalar polynomials of grade $k-1$. Then,
$\mathcal{B}_{n,d} = \{
n(z)^{k-1},
n(z)^{k-2}d(z),
\ldots ,
n(z)d(z)^{k-2} ,
d(z)^{k-1} \}$ is a basis of $\F[z]_{k-1}$. Therefore, there exists an invertible change of basis matrix $B\in \F^{k \times k}$ such that
\begin{equation} \label{eq:changeMobiusbasis} B \begin{bmatrix}
z^{k-1}\\
z^{k-2}\\
\vdots\\
z\\
1
\end{bmatrix} = \begin{bmatrix}
n(z)^{k-1}\\
n(z)^{k-2}d(z)\\
\vdots\\
n(z)d(z)^{k-2}\\
d(z)^{k-1}
\end{bmatrix} .
%\Leftrightarrow \mathcal{M}_{k-1} ( \sum_{i=0}^{k-1}v_i z^i) =  d(z)^{k-1} %\sum_{i=0}^{k-1} v_i \left( \frac{n(z)}{d(z)} \right)^i.
\end{equation}
\end{lemma}
\begin{proof} Let $r(z)=n(z)/d(z)$. The statement is a corollary of \cite[Theorem 3.18]{MacMMM14}, that says that  $p(z) \mapsto \mathcal{M}_{k-1,r(z)}(p)$ is an $\F$-linear automorphism of $\F[z]_{k-1}$.

%It is known  \cite[Theorem 3.18]{MacMMM14} that $p(z) \mapsto \mathcal{M}_{k-1,r(z)}(p)$ induces an $\F$-linear bijection, $\mathcal{M}_{k-1,r(z)} : \F[z]_{k-1} \longrightarrow \F[z]_{k-1}$. Thus, for any $q(z) \in \F[z]_{k-1}$, there exists exactly one $p(z) = \sum_{i=0}^{k-1} p_i z^i \in \F[z]_{k-1}$ such that $q(z) = \mathcal{M}_{k-1,r(z)}(p) = \sum_{i=0}^{k-1} p_i \, n(z)^i \, d(z)^{k-1-i}$. This implies that the set $\mathcal{B}_{n,d}$ with $k$ polynomials spans $\F[z]_{k-1}$. Hence, it is a basis because the dimension of $\F[z]_{k-1}$ is equal to $k$. Since $\{ z^{k-1}, z^{k-2}, \ldots , z , 1 \}$ is another basis of $\F[z]_{k-1}$, the existence of $B$ follows.
\end{proof}

\begin{theorem}\label{thm:infinitoacasa}
Let $P(z)\in \F[z]^{m \times n}$ be a matrix polynomial of grade $k$ and $v(z)$ be an ansatz polynomial of grade $k-1$. Let $n(z),d(z)$ be coprime linear scalar polynomials and denote by $\mathcal{M}_{g,r(z)} (X)$ the M\"{o}bius transform of grade $g$ of the matrix polynomial $X(z)$ induced by $r(z)=n(z)/d(z)$, that is, $\mathcal{M}_{g,r(z)}(X) = d(z)^g X(n(z)/d(z))$. Denote by $L(z)=z L_1+L_0$ the $\mathbb{DL} (P)$ pencil of $P(z)$ with ansatz polynomial $v(z)$, i.e., $\mathbb{DL} (P,v)$, and by $M(z)= z M_1 + M_0$ the $\mathbb{DL}(Q)$ pencil of $Q(z):=\mathcal{M}_{k,r(z)}(P)$  with ansatz polynomial $u(z):=\mathcal{M}_{k-1,r(z)}(v)$, i.e., $\mathbb{DL} (Q,u)$. Moreover, let $B$ be the $k \times k$ change of basis matrix in Lemma \ref{lemm:changebasemob}.
Then, $$(B^T \otimes I_m) \, \mathcal{M}_{1,r(z)} (L) \, (B \otimes I_n)=M(z).$$
In other words, the following diagram commutes:
$$
\begin{tikzpicture}[scale=2]
\node (A) at (0,2) {$(P,v)$};
\node (B) at (3,2) {$z L_1+L_0$};
\node (D) at (0,0) {$(Q,u)$};
\node (E) at (3,0) {$z M_1 + M_0$};
%\path[->,font=\scriptsize,>=angle 90]
\path[->,font=\scriptsize,]
(A) edge[bend right] node[right]{(k,k-1)-M\"{o}bius} (D)
(B) edge[bend left] node[left]{1-M\"{o}bius
+ $(B \otimes I)$-congruence} (E)
(A) edge node[above]{$\mathbb{DL}$} (B)
(D) edge node[above]{$\mathbb{DL}$} (E);
\end{tikzpicture}
$$
\end{theorem}
\begin{proof}
It is convenient, for the sake of clarity in the proof, to distinguish notationally between the variables associated with the triple $(P, v, L)$ and those associated with $(Q,u,M)$: to this goal, we will put hats on the variables associated with $P, v$ and $L$. Under the bijection \eqref{eq:bivariatebijection} and taking into account \eqref{eq:bezout},
$$L(\hat z) \mapsto \frac{P(\hat y)(\hat x-\hat z)v(\hat x)-P(\hat x)(\hat y-\hat z)v(\hat y)}{\hat x- \hat y}.$$
Taking into account \eqref{eq:bivariatebijection}, \eqref{eq:changeMobiusbasis}, and the fact that $\mathcal{M}_{k-1,r(z)}(\sum_{i=0}^{k-1} w_i z^i) = \sum_{i=0}^{k-1} w_i \, n(z)^i \, d(z)^{k-1-i}$, we see that left multiplication of $L(\hat z)$ by $B^T \otimes I_m$ corresponds functionally to a M\"{o}bius transformation of grade $k-1$ induced by $\hat y=r(y)$; similarly, right multiplication of $L(\hat z)$ by $B \otimes I_n$ corresponds to a M\"{o}bius transformation of grade $k-1$ induced by $\hat x=r(x)$; and finally the action of $\mathcal{M}_{1,r(z)}$ implements at the functional level a M\"{o}bius transformation of grade $1$ induced by $\hat z=r(z)$. Hence, using \eqref{eq:bezout}, the pencil $(B^T \otimes I) \mathcal{M}_{1,r(z)} (L) (B \otimes I)$ corresponds functionally to
\begin{equation}\label{lhs}
d(z) d(x)^{k-1} d(y)^{k-1}\frac{P(r(y))(r(x)-r(z))v(r(x))-P(r(x))(r(y)-r(z))v(r(y))}{r(x)-r(y)}.
\end{equation}
On the other hand, once again using \eqref{eq:bezout},
$$M(z) \mapsto \frac{Q(y)(x- z)u(x) -Q( x)(y- z)u(y)}{x - y},$$
or, equivalently after using the definitions of $Q$ and $u$, $M(z)$ is mapped by \eqref{eq:bivariatebijection} to the trivariate matrix polynomial
\begin{equation}\label{rhs}
\frac{d(y)^k d(x)^{k-1}P(r(y))(x- z)v(r(x)) -d(x)^kd(y)^{k-1}P(r(x))(y- z)v(r(y))}{ x - y}.
\end{equation}
Bearing in mind that \eqref{eq:bivariatebijection} is a bijection, proving the statement is therefore equivalent to showing that \eqref{lhs} and \eqref{rhs} are equal as trivariate polynomial functions in $x,y,z$. Dividing \eqref{lhs} and \eqref{rhs} by $d(y)^{k-1} d(x)^{k-1}$, the statement is then implied by
\begin{eqnarray*}
 \frac{P(r(y)) (d(z) r(x)-n(z))v(r(x))-P(r(x)) (d(z) r(y)-n(z))v(r(y))}{r(x)-r(y)} = \\
=\frac{d(y)P(r(y))(x-z)v(r(x))-d(x)P(r(x))(y-z)v(r(y))}{x-y}.
\end{eqnarray*}
Breaking each side as the difference of two similar terms, equating each pair of terms, and cancelling common factors, we see that in turn the latter equality is implied by
\begin{equation}\label{tres}
\begin{cases}
(x-y) (d(z)r(x)-n(z)) = (r(x)-r(y)) d(y)(x-z);  \\
(x-y)(d(z)r(y)-n(z)) = (r(x)-r(y)) d(x)(y-z).
\end{cases}
\end{equation}
We prove the first functional equation of the system \eqref{tres}, as the second one follows from the first one with the changes $x \mapsto y$ and $y \mapsto x$ and changing the signs of both sides. Multiplying both sides of the first equation in \eqref{tres} by $d(x)$ we obtain
$$(x-y) (d(z)n(x)-n(z)d(x)) = (n(x)d(y)-n(y)d(x)) (x-z),$$
whose correctness can be checked by using $n(s)=\alpha s + \beta , d(s)=\gamma s+ \delta$, and then expanding both sides to see that both are equal to
$(x-y)(x-z) (\delta \alpha - \beta \gamma)$.
%If either $x=y$ or $x=z$ (or both), the latter equality collapses to the %tautological statement that $0=0$; it therefore suffices to prove the equality %assuming that $x\ne y$ and $x\ne z$.
%To this goal, denote the bivariate Bezoutian function associated with $n$ and $d$ %by
%$$\beta(s,t)=\frac{n(t)d(s)-n(s)d(t)}{s-t};$$
%then it is clear that the sought functional equality is equivalent to
%$$\beta(z,x)=\beta(x,y).$$
%On the other hand, since by assumption $n$ and $d$ are coprime linear %polynomials, a straightforward calculation shows that $\beta(s,t)$ is a nonzero %constant. This concludes the proof.
\end{proof}

A consequence of Theorem \ref{thm:infinitoacasa} is that we can assume in the rest of the paper without loss of generality that both the matrix polynomial $P(z)$ and the ansatz polynomial $v(z)$ do not have eigenvalues and roots at infinity, which simplifies considerably the developments. Indeed, if $\infty$ is an eigenvalue of $P(z)$ or a root of $v(z)$, then we can apply a M\"{o}bius transformation induced by an appropriate $r(z)$ to both so that neither $Q(z)=\mathcal{M}_{k,r(z)} (P)$ nor $u(z)=\mathcal{M}_{k-1,r(z)}(v)$ has infinite eigenvalues/roots. Note that it is always possible to find such $r(z)$ because the set of eigenvalues/roots of $P$ and $v$ has a finite number of elements and $\F$ is algebraically closed, hence infinite. Taking into account Theorem \ref{thm:eigenstrucutreafterMobius}  and the fact that the strict equivalences $B^T \otimes I_m$ and $B \otimes I_n$ preserve the eigenvalues, the partial multiplicities and the minimal indices, Theorem \ref{thm:infinitoacasa} then allows us to obtain statements on the triple $(P,v,L=\mathbb{DL}(P,v))$ from those on the triple $(Q,u,M=\mathbb{DL}(Q,u))$.

\section{The minimal indices of $\mathbb{DL}(P,v)$} \label{sec-minimalbasisdlp}
In this section we will prove items 1 and 2 of Theorem \ref{thm:main}. We will prove first item 1 and then obtain item 2 as a simple corollary. The proof of item 1 proceeds by building a particular right minimal basis of the pencil $\mathbb{DL}(P,v)$ in \eqref{eq:bezout}. Throughout this section, we suppose that the, possibly singular, nonzero matrix polynomial $P(z) \in \F[z]^{m\times n}$ has grade $k\geq 2$, that the scalar ansatz polynomial $v(z)$ has grade $k-1$, and that the sets of eigenvalues of $P(z)$ and of roots of $v(z)$ are disjoint. We denote $\mathbb{DL}(P,v)=: z L_1 + L_0=:L(z)$. Moreover, having in mind to exploit Theorem \ref{thm:infinitoacasa}, we further suppose that %$P(z)$ has no infinite eigenvalues and that 
$v(z)$ has no infinite roots, which implies that %$P(z)$ has degree exactly $k$ and 
$v(z)$ has degree exactly $k-1$.% {\em For brevity, we will not repeat all these hypotheses in most of the results of this section.}

We start by reviewing some basic concepts on minimal bases that are used in this section and that can be found at \cite{For75,Mackey-filtration}. Recall that we are working in an algebraically closed field $\F$. For brevity, we say that a polynomial matrix $A(z) \in \F [z]^{n \times p}$, with $n\geq p$, is a minimal basis if its columns form a minimal basis of the rational subspace they span, which is denoted by $\spann A(z)$. Given a polynomial matrix $P (z) \in \F [z]^{m\times n}$, we say that $A(z)$ is a (right) minimal basis of $P(z)$ if $A(z)$ is a minimal basis and $\spann A(z) = \ker P(z)$. Let us write $A(z)$ in terms of its columns as  $A (z) = \begin{bmatrix} A_1(z) & \cdots &  A_p(z)\end{bmatrix}$. Then, we define the column-wise reversal of $A(z)$ as $(\cwRev A) (z) = \begin{bmatrix} \Rev A_1(z) & \cdots & \Rev A_p(z)\end{bmatrix}$, where we emphasize that each reversal is taken with respect to the degree of the corresponding column, and by convention $\deg 0=-\infty$. We denote $\{A(z)\}=(\cwRev A) (0)$. This constant matrix is termed in \cite{For75} the high order coefficient matrix of $A(z)$. The rules in the following lemma follow immediately from the definition of the high order coefficient matrix.

\begin{lemma}\label{lemm:columnwiserev} Let $A(z) \in \F [z]^{n \times p}$, $B(z) \in \F [z]^{n \times t}$ and $C \in \F^{n \times s}$, i.e., $C$ is a constant matrix. Then:
\begin{enumerate}
\item $\{\begin{bmatrix}A(z) & B(z)\end{bmatrix}\}=\begin{bmatrix}\{A(z)\} & \{B(z)\}\end{bmatrix}$;
\item $\{C\}=C$;
\item for any scalar polynomial $p(z) = p_d z^d + \cdots + p_1 z + p_0$, with $p_d \ne 0$,
\begin{enumerate}
\item $\{p(z) A(z)\}=p_d \, \{A(z)\}$;
\item if $\deg p(z) > \deg B(z)$ and no column of $A(z)$ is zero, then $\{p(z)A(z)+B(z)\}= p_d \, \{A(z)\}$.
\end{enumerate}
\end{enumerate}
\end{lemma}

Let us recall a very useful characterization of minimal bases that was first proved in \cite{For75}.

\begin{lemma}\label{forney}
Let $A(z) \in \mathbb{F}[z]^{n \times p}$, $n \geq p$. Then $A(z)$ is a minimal basis if and only if
\begin{enumerate}
\item $\rank A(z_0)=p$ for all $z_0 \in \mathbb{F}$ and
\item $\rank \{A(z)\}=p$.
\end{enumerate}
\end{lemma}

The following simple lemma will be also used in this section.

\begin{lemma} \label{lemm:blockrank}
Let $A,B,C$ be matrices of adequate sizes with entries in any field. Then
$\rank \begin{bmatrix}
0 & C \\
A & B
\end{bmatrix} = \rank \begin{bmatrix}
A & B\\
0 & C
\end{bmatrix} \geq \rank A + \rank C$.
\end{lemma}

In the following proposition, we give a first step towards constructing a right minimal basis of $\mathbb{DL}(P,v)$. Note also that, as a corollary of Proposition \ref{thetrick}, we obtain an alternative proof to that in \cite{DDM09} for the fact that none of the pencils in $\mathbb{DL}(P)$ is a linearization of $P(z)$ if $P(z)$ is singular, since the dimensions of the right kernels of $L(z)$ and $P(z)$ are different.

\begin{proposition}\label{thetrick}
Let $M(z) \in \F[z]^{n \times p}$ be a right minimal basis for $P(z)\in \F[z]^{m \times n}$ and suppose the distinct roots of $v(z)$ are $\mu_1,\dots,\mu_s \in \F$ with multiplicities $\ell_1,\dots,\ell_s$, respectively. Let $L(z) = \mathbb{DL}(P,v)$. Using the notation of Lemmata \ref{important}, \ref{confvanderpartone} and \ref{confvanderparttwo}, and taking into account that $\sum_{i=1}^{s} \ell_i = k-1$, define
\begin{align*} C & =\begin{bmatrix}
M_{\mu_1,k,\ell_1,M(z)} & M_{\mu_2,k,\ell_2,M(z)} & \cdots & M_{\mu_s,k,\ell_s,M(z)}
\end{bmatrix} \in \F^{kn \times (k-1)p}, \quad \mbox{and} \\
G (z) & = \begin{bmatrix} C & D(z) \end{bmatrix} \in \F[z]^{kn \times kp},
\end{align*}
where $D(z) := V(z) \otimes M(z) \in \F[z]^{kn \times p}$.
Then,
\begin{enumerate}
\item $C$ has full column rank and satisfies $L(z) \, C = 0$, that is $\spann C \subseteq \ker L(z)$. Moreover, $L(z)$ has at least $(k-1)p$ right minimal indices equal to $0$,  where $p = \dim \ker P(z)$;
\item $G(z)$ has full column rank and satisfies $L(z) \, G(z) = 0$, that is $\spann G(z) \subseteq \ker L(z)$.
\end{enumerate}
\end{proposition}

\begin{proof} Item 1. Taking into account Lemmata \ref{confvanderpartone} and \ref{confvanderparttwo}, we consider $C$ partitioned as
$C = \begin{bmatrix} C_1 & C_2 & \cdots & C_{k-1} \end{bmatrix}$, where each $C_i$ has a structure of the type $ C_i = \left. [V(z) \otimes M (z)]^{(a)}\right|_{z=\mu_j} \in \F^{kn \times p}$ for some root $\mu_j$ of $v(z)$ and some $a$ such that $0\leq a \leq  \ell_j -1$. Then $L(z)C= \begin{bmatrix} L(z) C_1 & \cdots & L(z) C_{k-1} \end{bmatrix}$. Since, according to \eqref{eq:bezout}, $L(z) \mapsto \mathcal{B}_{P,v} (x,y,z)$ via the bivariate matrix polynomial bijection \eqref{eq:bivariatebijection}, right multiplication of $L(z)$ by the block column $C_i$ corresponds in the functional bivariate formulation to right multiplication of $\mathcal{B}_{P,v} (x,y,z)$ by $M(x)$, followed by taking the $a$-th partial derivative with respect to $x$ and $x$-evaluation at $x=\mu_j$ (see also \cite[Table 3.1]{NNT17}). Explicitly, this means
 $$L(z) C_i \mapsto \left. \frac{\partial^a}{\partial x^a} \frac{P(y)(x-z)v(x) - (y-z)v(y) P(x)}{x - y} M(x) \right|_{x=\mu_j}.$$
Since $P(x) M(x) = 0$, we only need to pay attention to the derivative of the first term and, taking into account that $\mu_j$ is a root of $v(x)$ of multiplicity $\ell_j > a$ (crucially, the inequality is strict), we then see that $L(z) C_i \mapsto 0$ for all $i=1,\ldots k-1$, and hence $L(z) C=0$ because \eqref{eq:bivariatebijection} is a bijection\footnote{We warn the reader that, as explained in \cite[Section 3]{NNT17}, the bivariate polynomial bijection \eqref{eq:bivariatebijection} can clearly be defined more in general from $k\times h$ block matrices with blocks of size $m\times n$ to $m\times n$ bivariate matrix polynomials of grade $k-1$ in $y$ and of grade $h-1$ in $x$.}.

To prove that $C$ has full column rank, it suffices to use Lemma \ref{forney}-1 on $M(z)$ and then to invoke (a special case of) Lemma \ref{confvanderparttwo}. Finally, since the columns of $C$ are linearly independent constant vectors and belong to $\ker L(z)$, they can be completed to a minimal basis of $\ker L(z)$ via the method described in \cite[Ch. XII, \S 5]{gantmacher}. Thus, the $(k-1)p$ columns of the constant matrix $C$ give $(k-1)p$ minimal indices of $L(z)$ equal to zero. Another way to get the same result is through the so-called ``Strong Minimality Property of Minimal Indices'' proved in \cite[Theorem 4.2]{Mackey-filtration}.

Item 2. Observe that the operation $L(z) D(z)$ in the bivariate matrix polynomial framework amounts to multiply $\mathcal{B}_{P,v} (x,y,z)$ by $M(x)$ on the right followed by $x$-evaluation at $x = z$. Explicitly, this means
$$L(z)D(z) \mapsto \frac{P(y)(z-z)v(z) - (y-z)v(y) P(z)}{z - y} M(z) = 0,$$
where we have used again that $P(z)M(z) = 0$. Therefore, $L(z) D(z) = 0$, which combined with the result in item 1 implies $L(z) G(z) =0$. To see that $G(z)$ has full column rank, we apply Lemma \ref{confvanderparttwo} to $G(\alpha)$, where $\alpha$ is any element of $\F$ such that $\alpha \neq \mu_i$, for $i=1, \ldots , s$. This proves that the constant matrix $G(\alpha)$ has full column rank, which implies that the matrix polynomial $G(z)$ has full column rank.
\end{proof}

Proposition \ref{thetrick} implies that $\dim \ker L(z) \geq kp$. In fact, we will prove later that $\dim \ker L(z) = kp$. This, and the fact that $G(z)$ is not a minimal basis, motivates us to look for a minimal basis $F(z) \in \F [z]^{kn \times kp}$ such that $\spann F(z) = \spann G(z)$. Observe that $G(z)$ is not a minimal basis because it does not satisfy the first condition of Lemma \ref{forney}: indeed, $D(\mu_i)$ is equal to the first block column of $M_{\mu_i,k,\ell_i,M(z)}$ for $i=1,\ldots, s$. A minimal basis $F(z)$ is described in  Proposition \ref{prop.minimalF} in terms of the Hermite interpolating polynomials (see, for instance, \cite[p. 53]{stoer-bulirsch}) associated with the distinct roots $\mu_1 , \ldots , \mu_s$ of $v(z)$ and their multiplicities $\ell_1, \ldots , \ell_s$.

\begin{proposition} \label{prop.minimalF} Let $M(z) \in \F[z]^{n \times p}$ be a right minimal basis for $P(z)\in \F[z]^{m \times n}$ and let $\gamma_1 \leq \cdots \leq \gamma_p$ be the degrees of the columns of $M(z)$, i.e., $\gamma_1 \leq \cdots \leq \gamma_p$ are the right minimal indices of $P(z)$. Suppose the distinct roots of $v(z)$ are $\mu_1,\dots,\mu_s \in \F$ with multiplicities $\ell_1,\dots,\ell_s$, respectively. Let $L(z) = \mathbb{DL}(P,v)$ and let $C, D(z)$ and $G(z)$ be the matrices defined in Proposition \ref{thetrick}. Let us partition the matrix $C \in \F^{kn \times (k-1)p}$ as $C = \begin{bmatrix}
C_1 & C_2 & \cdots & C_{k-1} \end{bmatrix}$ with $C_i \in \F^{kn \times p}$ for $i = 1 , \ldots , k-1$. Let $H_1 (z), H_2 (z) , \ldots , H_{k-1} (z)$ be the scalar Hermite interpolating polynomials associated with $\mu_1 , \ldots , \mu_s$ and $\ell_1, \ldots , \ell_s$ ordered as follows: if $ C_i = \left. [V(z) \otimes M (z)]^{(a)}\right|_{z=\mu_j} \in \F^{kn \times p}$ for some $\mu_j$ and some $a$ such that $0\leq a \leq  \ell_j -1$, then
\begin{equation} \label{eq.hermitepolys}
H_i^{(b)} (\mu_t)  = \begin{cases}
1, \ & \ \mathrm{if} \; t = j \; \mbox{and} \; b = a,\\
0,  \ & \ \mathrm{if} \; t \neq j \; \mbox{or} \; b \ne a,
\end{cases}
\end{equation}
where, for each $\mu_t$,  the integer $b$ satisfies $0 \leq b \leq \ell_t - 1$. Define
\[ E(z) = D(z) - \sum_{i=1}^{k-1} H_i(z)\, C_i \in \F[z]^{kn \times p} \quad \mbox{and} \quad
F (z)  = \begin{bmatrix} C & \frac{1}{v(z)} E(z) \end{bmatrix} .
\]
Then,
\begin{enumerate}
  \item $F(z)$ has full column rank, satisfies $L(z) \, F(z) = 0$, that is $\spann F(z) \subseteq \ker L(z)$, and $\spann F(z) = \spann G(z)$;
  \item $F(z)\in \F[z]^{kn \times kp}$, i.e., it is a matrix polynomial, and the degrees of its columns are $\underbrace{0,\dots,0}_{(k-1)p},\gamma_1,\dots,\gamma_p$; and
  \item $F(z)$ is a minimal basis.
\end{enumerate}
\end{proposition}
\begin{proof}
Item 1. Observe that $F(z) = G(z) \, S(z)$, where
\begin{equation} \label{eq.factorizationFz}
S(z) = \left[
\begin{array}{ccccc}
  I_p &  &  &  & -H_1(z) I_p \\
    & I_p &  &  & -H_2 (z) I_p \\
      &  & \ddots &  & \vdots \\
      &  &  & I_p & -H_{k-1}(z) I_p \\
      &  &  &   & I_p
\end{array} \right] \,
\left[
\begin{array}{ccccc}
  I_p &  &  &  &  \\
    & I_p &  &  & \\
      &  & \ddots &  & \\
      &  &  & I_p &  \\
      &  &  &   & \frac{1}{v(z)} I_p
\end{array} \right]
\end{equation}
is an invertible matrix over $\F(z)$. This relation and the properties of $G(z)$ established in item 2 of Proposition \ref{thetrick} immediately imply the result.

Item 2. We only need to prove that $\frac{1}{v(z)} E(z)$ is a matrix polynomial and to determine the degrees of its columns. From the definition of $D(z)$, it is obvious that the degrees of its columns are $\gamma_1 + k-1,\dots,\gamma_p + k-1$. On the other hand, the Hermite interpolating polynomials $H_1 (z), \ldots , H_{k-1} (z)$ have degrees at most $k-2$, since each of them is defined through $k-1$ interpolatory conditions. Thus, the columns of $E(z)$ have degrees $\gamma_1 + k-1,\dots,\gamma_p + k-1$. The definitions of $D(z), H_i(z)$ and $C_i$ imply that
\[
E^{(b)} (\mu_t) = D^{(b)}(\mu_t) - \sum_{i=1}^{k-1} H_i^{(b)} (\mu_t) \, C_i = 0, \quad \mbox{for} \quad t =1,\ldots  s, \; b = 0,\ldots , \ell_t - 1,
\]
which in turn implies that $v(z)$ divides each entry of $E(z)$. Therefore, $\frac{1}{v(z)} E(z)$ is a matrix polynomial and the degrees of its columns are $\gamma_1,\dots,\gamma_p$, since the degree of $v(z)$ is $k-1$.

Item 3. We will prove that $F(z)$ is a minimal basis by using Lemma \ref{forney}. We start by showing that $F(z)$ satisfies the second condition of Lemma \ref{forney}. For this purpose, denote $A(z) = \frac{1}{v(z)} E(z)$ and let $v_{k-1}$ be the (nonzero) coefficient of $z^{k-1}$ in $v(z)$. Then, Lemma \ref{lemm:columnwiserev}-3 implies $\{E(z)\} = \{ v(z) \, A(z) \} = v_{k-1} \{A(z)\}$. Using again that the polynomials $H_i (z)$ have degree at most $k-2$, we get that $\{E(z)\} = \{D(z)\} = \{V(z) \otimes M(z)\} = \left[
\begin{smallmatrix} \{M(z)\} \\ 0 \end{smallmatrix}\right]$, due to the structure of $V(z)$ (recall Lemma \ref{important}). Combining these results and using again Lemma \ref{lemm:columnwiserev}, we get
\begin{equation} \label{eq.minFzargument}
\{F(z)\}=\{\begin{bmatrix} C & \frac{1}{v(z)} E(z)\end{bmatrix}\} = \begin{bmatrix} \{C\} & \frac{1}{v_{k-1}}\{E(z)\}\end{bmatrix} = \begin{bmatrix} \tilde{C} & \frac{1}{v_{k-1}} \{M(z)\}\\
\hat{C} & 0\end{bmatrix},
\end{equation}
where $\tilde{C}$ (resp., $\hat{C}$) contains the first $n$ rows (resp., the rows from $(n+1)$-th to last) of $C$. Note that taking into account the definition of $C$,
$$\hat{C} = \begin{bmatrix}
M_{\mu_1,k-1,\ell_1,M(z)} & M_{\mu_2,k-1,\ell_2,M(z)} & \cdots & M_{\mu_s,k-1,\ell_s,M(z)}
\end{bmatrix} \in \F^{(k-1)n \times (k-1)p},$$ which
has full column rank by Lemma \ref{confvanderparttwo}, since $\sum_{i=1}^{s} \ell_i = k-1$. Moreover, $\{M(z)\}$ has full column rank because $M(z)$ is a minimal basis. Hence, \eqref{eq.minFzargument} implies that $\{F(z)\}$ has full column rank.

For proving that $F(z)$ satisfies the first condition of Lemma \ref{forney}, we distinguish two cases: $z_0$ is not a root of $v(z)$ and $z_0$ is a root of $v(z)$. Suppose first that $v(z_0) \neq 0$. Then $G(z_0) = \begin{bmatrix} C & D(z_0) \end{bmatrix}$ has full column rank by Lemma \ref{confvanderparttwo} (we have already seen this in the proof of item 2 of Proposition \ref{thetrick}). Moreover, using \eqref{eq.factorizationFz}, we get $F(z_0) = G(z_0) \, S(z_0)$, with $S(z_0)$ invertible. This implies that $F(z_0)$ and $G(z_0)$ have the same rank and, hence, that $F(z_0)$ has full column rank.

Next, we assume that $v(z_0) = 0$, i.e., $z_0 = \mu_i$ for some $i = 1, \ldots , s$. Without loss of generality, we can assume $z_0 = \mu_1$, since otherwise we can simply permute the columns of $C$ and place $M_{\mu_i,k,\ell_i,M(z)}$ as the first block column. In this case, we evaluate $E(z)/v(z)$ at $\mu_1$ by applying De L'Hospital rule $\ell_1$ times. This yields
\[ \left[ \frac{E(z)}{v(z)} \right]_{z=\mu_1} =   \frac{1}{v^{(\ell_1)}(\mu_1)} \,
\left(
D^{(\ell_1)}(\mu_1) - \sum_{i=1}^{k-1} H_i^{(\ell_1)}(\mu_1) \, C_i \right) , \]
where we note that $v^{(\ell_1)}(\mu_1) \neq 0$. Then,
$$
F(\mu_1) = \begin{bmatrix} C &  D^{(\ell_1)}(\mu_1) \end{bmatrix} \,
\begin{bmatrix} I_p & & & & -H^{(\ell_1)}_1(\mu_1)I_p  \\
 & I_p & & & -H^{(\ell_1)}_2(\mu_1) I_p\\
& & \ddots & & \vdots\\
& & & I_p & -H^{(\ell_1)}_{k-1}(\mu_1) I_p\\
& & & & I_p
\end{bmatrix} \left[
\begin{array}{ccccc}
  I_p &  &  &  &  \\
    & I_p &  &  & \\
      &  & \ddots &  & \\
      &  &  & I_p &  \\
      &  &  &   & \frac{1}{v^{(\ell_1)}(\mu_1)} I_p
\end{array} \right].
$$
The last two matrices in the equation above are invertible. Thus
$$\rank F(\mu_1) = \rank \begin{bmatrix} C &  D^{(\ell_1)}(\mu_1) \end{bmatrix} = \rank \begin{bmatrix} C &  \left. [ V(z) \otimes M(z)]^{(\ell_1)} \right|_{z = \mu_1} \end{bmatrix}.$$
Note that, taking into account the definition of $C$, we can permute the columns of
$\begin{bmatrix} C &  D^{(\ell_1)}(\mu_1) \end{bmatrix}$ for moving the block column $D^{(\ell_1)}(\mu_1)$ just after the block $M_{\mu_1,k,\ell_1,M(z)}$ of $C$ to get
$$
\begin{bmatrix}
M_{\mu_1,k,\ell_1 + 1,M(z)} & M_{\mu_2,k,\ell_2,M(z)} & \cdots & M_{\mu_s,k,\ell_s,M(z)}
\end{bmatrix} \in \F^{kn \times kp}.
$$
Since $\ell_1 + 1 + \sum_{i=2}^{s} \ell_i = k$, this matrix has full column rank by Lemma \ref{confvanderparttwo}, which concludes the proof.
\end{proof}

Observe that if $\dim \ker L(z) = kp$, then Proposition \ref{prop.minimalF} implies that $F(z)$ is a right minimal basis for $L(z)$, which proves item 1 in Theorem \ref{thm:main} taking into account the degrees of the columns of $F(z)$. Note that this proof would work in the case $v(z)$ has no infinite roots %and $P(z)$ has no infinite eigenvalues
, since we are assuming %these conditions 
this condition from the beginning of this section. Therefore, our next task is to prove that $\dim \ker L(z) = kp$, which will be a consequence of combining Theorem \ref{blockevaluation} with Proposition \ref{prop.minimalF}. This will lead to the main result of this section, that is, Theorem \ref{thm:Fisminimal}.

Before stating and proving Theorem \ref{blockevaluation}, we introduce the auxiliary Definition \ref{def.auxWconfluentVander} and the technical Lemma \ref{lem:combinatorial}.

\begin{definition} \label{def.auxWconfluentVander} Let $w(z)\in \F[z]$ be a scalar polynomial of degree $k$ whose distinct roots are $\mu_0, \mu_1 , \ldots , \mu_s \in \F$ with multiplicities $\ell_0, \ell_1 , \ldots , \ell_s$. Let $V(z) \in \F[z]^{k \times 1}$ be the Vandermonde vector of degree $k-1$ introduced in Lemma \ref{important}. Let us define the constant matrices
$$
W_{\mu_i} = \begin{bmatrix} V(\mu_i) & V^{(1)}(\mu_i) & \displaystyle  \frac{1}{2!} V^{(2)}(\mu_i) & \cdots & \displaystyle \frac{1}{(\ell_i -1)!} V^{(\ell_i - 1)}(\mu_i)  \end{bmatrix} \in \F^{k \times \ell_i} \quad \mbox{for $i=0, 1, \ldots, s$}.
$$
The normalized confluent Vandermonde matrix associated with the polynomial $w(z)$ is
\[
W = \begin{bmatrix}
      W_{\mu_0} & W_{\mu_1} & \cdots & W_{\mu_s}
    \end{bmatrix} \in \F^{k \times k}.
\]
\end{definition}

The normalized confluent Vandermonde matrix is obviously invertible, since it is just a standard confluent Vandermonde matrix associated with a set of distinct nodes with its columns multiplied by nonzero scalars. In fact, we might use the standard confluent Vandermonde matrix in our developments, but the normalization via the factorials will simplify the algebraic manipulations. In this context, here and below, the $j$th normalized derivative is defined as
$$ \hat \partial^j_{x^j} = \frac{1}{j!} \frac{\partial^j}{\partial x^j}.$$
\begin{lemma}\label{lem:combinatorial}
Let $\mu \in \F$ and let $h_c(x,y)$ be the $c$th complete homogenous polynomial,
\[ h_c(x,y) = \sum_{h=0}^c x^{c-h} y^h.\]
Then, $\hat \partial^a_{x^a} \hat \partial^b_{y^b} \, h_c(x,y)$ evaluated at $x=y=\mu$ only depends on $a+b$, and it is equal to
\[  \left[ \hat \partial^a_{x^a} \hat \partial^b_{y^b} h_c(x,y)  \right]_{x=y=\mu}  = \begin{cases}
0 \ & \ \mathrm{if} \ a+b > c;\\
\begin{pmatrix}
c+1\\
a+b+1
\end{pmatrix} \mu^{c-a-b} \ & \ \mathrm{if} \ a+b \leq c.
\end{cases}  \]
\end{lemma}
\begin{proof}
We may assume that $a+b \leq c$, as the statement is otherwise obvious. A direct computation yields
\[
 \hat \partial^a_{x^a} \hat \partial^b_{y^b} h_c(x,y) = \sum_{h=0}^{c}
 \begin{pmatrix}
   c-h \\
   a
 \end{pmatrix} \,
\begin{pmatrix}
   h \\
   b
 \end{pmatrix}  \, x^{c-h-a} \, y^{h-b},
\]
where $\left( \begin{smallmatrix} \ell \\ k \end{smallmatrix} \right) = 0$ if $\ell < k$, as usual. Thus
\[
\left[ \hat \partial^a_{x^a} \hat \partial^b_{y^b} h_c(x,y)  \right]_{x=y=\mu}  =
\mu^{c-a-b} \, \sum_{h=0}^{c}
 \begin{pmatrix}
   c-h \\
   a
 \end{pmatrix} \,
\begin{pmatrix}
   h \\
   b
 \end{pmatrix} = \begin{pmatrix}
c+1\\
a+b+1
\end{pmatrix} \mu^{c-a-b},
\]
where the last equality follows from \cite[eq. (4)]{adsokal}, which is a corollary of the famous Chu–Vandermonde identity.
\end{proof}

We are now ready to prove a result that will allow us to get a lower bound on $\rank L(z)$ or, equivalently, due to the rank-nullity theorem, an upper bound on $\dim \ker L(z)$.

\begin{theorem}\label{blockevaluation} Let $P(z)\in \F[z]^{m\times n}$ be a matrix polynomial of degree $k\geq 2$, $v(z)$ be a scalar ansatz polynomial of degree $k-1$, and $L(z) = \mathbb{DL}(P,v) \in \F [z]^{km \times kn}$. Suppose the distinct roots of $v(z)$ are $\mu_1, \ldots, \mu_s$ with multiplicities $\ell_1, \ldots, \ell_s$ and let $\mu_0 \in \F$ be any element such that $\mu_0 \ne \mu_i$ for $i = 1,\ldots , s$. Define the scalar polynomial $w(z) = (z-\mu_0) \, v(z)$ of degree $k$ and let $W\in \F^{k\times k}$ be the normalized confluent Vandermonde matrix associated with the polynomial $w(z)$, with $\ell_0 =1$. Then, there exist nonzero constants $c_0,c_1, \dots, c_s \in \F$ equal to
$$ c_i = \frac{1}{\ell_i!} \, w^{(\ell_i)} (\mu_i),$$
 such that
\[     (W^T \otimes I_m)    L(\mu_0)       (W \otimes I_n) = \bigoplus_{i=0}^s Q_i, \qquad Q_i = \begin{bmatrix}
0 & 0 & \dots & 0 & c_i P(\mu_i)\\
0 & 0 &  &c_i P(\mu_i) & \star\\
\vdots & \vdots & \iddots & & \\
0 & c_i P(\mu_i) & \star & \dots & \star\\
c_i P(\mu_i) & \star & \star & \dots & \star
\end{bmatrix}  \in \F^{m \ell_i \times n \ell_i} ,  \]
i.e., each $Q_i$ is block antitriangular with all its antidiagonal blocks equal to $c_i P(\mu_i)$. The symbol $\star$ denotes a block whose precise form is not specified.
\end{theorem}
\begin{proof} We consider the constant matrix $(W^T \otimes I_m)    L(\mu_0)       (W \otimes I_n) \in \F^{km \times kn}$ partitioned into $k\times k$ blocks of size $m\times n$. In order to determine these blocks recall that, according to \eqref{eq:bezout}, $L(\mu_0) \mapsto \mathcal{B}_{P,v} (x,y,\mu_0)$ via the bijection \eqref{eq:bivariatebijection}. Moreover, the $(\alpha,\beta)$ $m\times n$ block of $(W^T \otimes I_m) L(\mu_0) (W \otimes I_n)$ is just the $\alpha$th block row of $(W^T \otimes I_m)$ (of size $m \times km$) times $L(\mu_0)$ times the $\beta$th block column of $(W \otimes I_n)$ (of size $kn \times n$). But taking into account the structure of the columns of $W$ and the equation \eqref{eq:bivariatebijection}, left multiplication of $L(\mu_0)$ by any block row of $(W^T \otimes I_m)$ corresponds to perform on $\mathcal{B}_{P,v} (x,y,\mu_0)$ a $b$th normalized derivative with respect to $y$, $0 \leq b < \ell_j$, followed by evaluation at some $y=\mu_j$. Similarly, right multiplication of $L(\mu_0)$ by any block column of $(W \otimes I_n)$ corresponds to perform on $\mathcal{B}_{P,v} (x,y,\mu_0)$ an $a$th normalized derivative with respect to $x$, $0 \leq a < \ell_i$, followed by evaluation at some $x=\mu_i$. Hence, each $m\times n$ block of $(W^T \otimes I_m)    L(\mu_0)       (W \otimes I_n) \in \F^{km \times kn}$ has the form
\begin{equation} \label{eq.blocksrevealingrank}
\left[   \hat \partial^a_{x^a} \hat \partial^b_{y^b}   \frac{P(y)w(x)-w(y)P(x)}{x-y}    \right]_{x=\mu_i,y=\mu_j},
\end{equation}
for $0 \leq j,i \leq s$, $0\leq b \leq \ell_j - 1$ and  $0\leq a \leq \ell_i - 1$.
Suppose first $i \neq j$ so that $\mu_i \neq \mu_j$. Then, for all fixed $y \neq \mu_i$, the function $(x-y)^{-1}P(y)w(x)$ has a zero in $x$ at $\mu_i$ of multiplicity $\ell_i$, and hence its $a$th normalized derivative with respect to $x$ is zero when evaluated at $x=\mu_i$. Analogously, the function $(x-y)^{-1} w(y) P(x)$, for all fixed $x \neq \mu_j$, has a zero in $y$ at $\mu_j$ of multiplicity $\ell_j$, and therefore its $b$th normalized derivative with respect to $y$ is zero when evaluated at $y=\mu_j$. We conclude that the corresponding block is $0$. Hence
\[
(W^T \otimes I_m) L(\mu_0) (W \otimes I_n) = \bigoplus_{i=0}^s Q_i, \qquad Q_i =
(W_{\mu_i}^T \otimes I_m) L(\mu_0) (W_{\mu_i} \otimes I_n),
\]
where we have used the notation in Definition \ref{def.auxWconfluentVander}.

The case $i=j$ corresponds to analyzing the blocks within $Q_i$ with  \eqref{eq.blocksrevealingrank} evaluated at $x=y=\mu_i$. Each of these blocks corresponds to a pair $(b,a)$ with $0 \leq a,b \leq \ell_i-1$, where $\ell_i$ is the multiplicity of $\mu_i$ as a root of $w(z)$. Note that we are only interested in the case where $a+b \leq \ell_i -1$, since the blocks corresponding to $a+b > \ell_i -1$ are not specified in the statement. Recall that $b$ corresponds to the row-block index and $a$ to the column-block index.

Define
\[ \Pi(x,y)=\frac{P(x)-P(y)}{x-y}, \qquad \Omega(x,y)=\frac{w(x)-w(y)}{x-y},\]
so that
\begin{equation} \label{eq.twoaddends}
\frac{P(y)w(x)-w(y)P(x)}{x-y} = P(y) \Omega(x,y) - w(y) \Pi(x,y).
\end{equation}
Since differentiation is linear, we can consider the two addends in the right-hand side of \eqref{eq.twoaddends} separately. Note first that
\[ \hat \partial^a_{x^a} \hat\partial^b_{y^b} [w(y) \Pi(x,y)] = \frac{1}{a!b!} \sum_{h=0}^b \begin{pmatrix}
b\\
h
\end{pmatrix} w^{(h)}(y)\,  \partial^a_{x^a} \cdot \partial^{b-h}_{y^{b-h}} \Pi(x,y);\]
since we eventually evaluate at $x=y=\mu_i$, and since $w^{(h)}(\mu_i)=0$ for all $0\leq h \leq b < \ell_i$, the contribution of the second addend to any of the $m\times n$ blocks we are investigating is $0$. For the first addend in the right-hand side of \eqref{eq.twoaddends}, we observe that
\[ \hat \partial^a_{x^a} \hat\partial^b_{y^b} [P(y) \Omega(x,y)] = \sum_{h=0}^{b} \frac{1}{(b-h)!} P^{(b-h)}(y) \, \hat \partial^a_{x^a} \hat \partial^{h}_{y^h} \Omega(x,y). \]
On the other hand, using the relation $x^{c+1}-y^{c+1}=(x-y)h_c(x,y)$, where $h_c(x,y)$ is the $c$th complete homegenous polynomial in Lemma \ref{lem:combinatorial}, we see that
$$ \Omega(x,y) = \sum_{c=0}^{k-1} w_{c+1} h_c(x,y),$$
where $w(z) = \sum_{t=0}^{k} w_t \, z^t$.
Hence, by Lemma \ref{lem:combinatorial},
$$[\hat \partial^a_{x^a} \hat \partial^{h}_{y^h} \Omega(x,y)]_{x=y=\mu_i} = \sum_{c=a+h}^{k-1} \begin{pmatrix}
c+1\\
a+h+1
\end{pmatrix} w_{c+1} \mu_i^{c-a-h} = \frac{1}{(a+h+1)!} w^{(a+h+1)} (\mu_i).$$
It follows that the $(b,a)$ block of $Q_i$ has the expression
\begin{equation}\label{eq:dablock}
\left[   \hat \partial^a_{x^a} \hat \partial^b_{y^b}   \frac{P(y)w(x)-w(y)P(x)}{x-y}    \right]_{x=y=\mu_i} =
 \sum_{h=0}^b \frac{1}{(b-h)!} P^{(b-h)} (\mu_i) \frac{1}{(a+h+1)!} w^{(a+h+1)} (\mu_i).
\end{equation}
If $a+b < \ell_i -1$, then $a+h+1 \leq a+b+1 < \ell_i$ and hence $w^{(a+h+1)} (\mu_i) = 0$. This shows that all the blocks above the main antidiagonal are zero. For the blocks on the antidiagonal, the only possibility for the addend in the outer summation to be nonzero is that $h=b$ so that $a+h+1=a+b+1=\ell_i$. Plugging this information into \eqref{eq:dablock}, and taking into account the definition of $c_i$, we obtain $ c_i P(\mu_i)$ as desired. To conclude the proof we need to show that $c_i \neq 0$. This must be true, for if not the multiplicity of $\mu_i$ as a root of $w$ would not be $\ell_i$.
\end{proof}

As announced, we are now in the position to state and prove the main result of this section. %For completeness, we include all the necessary hypotheses mentioned in the first paragraph of Section \ref{sec-minimalbasisdlp}.

\begin{theorem} \label{thm:Fisminimal} Let $P(z) \in \F [z]^{m\times n}$ be a matrix polynomial of grade $k \geq 2$ and $v(z)$ be a scalar ansatz polynomial of grade $k-1$. Suppose that the sets of eigenvalues of $P(z)$ and of roots of $v(z)$ are disjoint, %that $P(z)$ has no infinite eigenvalues 
and that $v(z)$ has no infinite roots. Let $L(z) = \mathbb{DL} (P,v) \in \F [z]^{km \times kn}$, $p = \dim \ker P(z)$ and $F(z) \in \F [z]^{kn \times kp}$ be the matrix polynomial in Proposition \ref{prop.minimalF}. Then $F(z)$ is a minimal basis for $\ker L(z)$, i.e., $F(z)$ is a right minimal basis of $L(z)$.
\end{theorem}
\begin{proof} In Proposition \ref{prop.minimalF}, we have already proved that $L(z) F(z)=0$ (or, equivalently, that $\spann F(z) \subseteq \ker L(z)$), that $F(z)$ has full column rank and that $F(z)$ is a minimal basis. Therefore, it only remains to prove that the dimension of $\ker L(z)$ is precisely $kp$, i.e., that it is equal to the rank of $F(z)$.

From the results above, we obviously have that $\dim \ker L(z) \geq kp = \dim \spann F(z)$. We now argue that $\rank L(z) \geq k(n-p)$, which is equivalent to $\dim \ker L(z) \leq kp$ by the rank-nullity theorem, and therefore implies $\dim \ker L(z) = kp$. To this goal, we use Theorem \ref{blockevaluation}, choosing an element $\mu_0$ which is not an eigenvalue of $P(z)$, and recall that the roots $\mu_i$, $i =1, \ldots, s$, of $v(z)$ are not eigenvalues of $P(z)$ by the assumption that the sets of eigenvalues of $P(z)$ and of roots of $v(z)$ are disjoint. Thus, $\rank P(\mu_i) = \rank P(z) = n-p$, for $i =0,1 ,\ldots , s$. Then, by Theorem \ref{blockevaluation}, by Lemma \ref{lemm:blockrank}, and by the fact that the matrix $W$ in Theorem \ref{blockevaluation} is nonsingular, we have, using the notation of Theorem \ref{blockevaluation},
 $$ \rank L(z) \geq \rank L(\mu_0) = \sum_{i=0}^s \rank Q_i \geq \sum_{i=0}^s \ell_i \rank P(\mu_i) = (n-p) \sum_{i=0}^s \ell_i = (n-p)k.$$
 \end{proof}

As a corollary of Theorem \ref{thm:Fisminimal} and the properties of $\mathbb{DL}(P)$ pencils, we prove items 1 and 2 of Theorem \ref{thm:main}.

\begin{theorem}\label{thm:Parts1and2main} {\rm (Items 1 and 2 of Theorem \ref{thm:main})} Let $P(z) \in \F[z]^{m \times n}$ be a matrix polynomial of grade $k \geq 2$, with $\dim \ker P(z) = p$ and $\dim \ker P(z)^T = q$, and let $v(z) \in \F[z]$ be a scalar polynomial of grade $k-1$. Suppose that the sets of eigenvalues of $P(z)$ and of roots of $v(z)$, each including possibly infinite eigenvalues or roots, are disjoint. Denote by $L(z)\in \F [z]^{km \times kn}$ the pencil in $\mathbb{DL}(P)$ with ansatz polynomial $v(z)$, i.e., $L(z) = \mathbb{DL} (P,v)$. Then
\begin{enumerate}
\item The right minimal indices of $P(z)$ are $\gamma_1 \leq \dots \leq \gamma_p$ if and only if the right minimal indices of $L(z)$ are $\underbrace{0=\dots=0}_{p(k-1) \ \mathrm{times}} \leq \gamma_1 \leq \dots \leq \gamma_p$;
\item The left minimal indices of $P(z)$ are $\eta_1 \leq \dots \leq \eta_q$ if and only if the left minimal indices of $L(z)$ are $\underbrace{0=\dots=0}_{q(k-1) \ \mathrm{times}} \leq \eta_1 \leq \dots \leq \eta_q$.
\end{enumerate}
\end{theorem}
\begin{proof}
Item 1. Taking into account Theorems \ref{thm:eigenstrucutreafterMobius} and \ref{thm:infinitoacasa}, we can assume without loss of generality %that $P(z)$ has no infinite eigenvalues and
 that $v(z)$ has no infinite roots, since, otherwise, we can apply a M\"{o}bius transformation to $P(z)$ and $v(z)$ to get a matrix polynomial $Q(z)$ and an ansatz polynomial $u(z)$ without %eigenvalues and
  roots at infinity%, respectively
  , and work with $Q(z)$ and the pencil $\mathbb{DL}(Q,u)$ without changing the minimal indices of $P(z)$ and of $\mathbb{DL} (P,v)$ (recall the last paragraph in Subsection \ref{subsec-mobiusdlp}). Then, the result follows from Theorem \ref{thm:Fisminimal} and item 2 in Proposition \ref{prop.minimalF}.

Item 2. The left minimal indices of $P(z)$ and of $L(z)$ are the right ones of $P(z)^T$ and of $L(z)^T$, respectively. With this in mind, the result follows immediately from the result in item 1 and the equality $L(z)^T = ( \mathbb{DL}(P,v))^T = \mathbb{DL}(P^T,v)$. When $P(z)$ is square, this is a consequence of the fact that the pencils in $\mathbb{DL}(P)$ are block symmetric \cite{hmmtsymlin,MacMMM06,NNT17}. For more general rectangular matrix polynomials, it can be proved as follows. Let $B$ be any $k\times k$ block matrix with blocks $B_{ij} \in \F^{m\times n}$ and $F(x,y)$ be its associated $m\times n$ bivariate matrix polynomial via the bijection \eqref{eq:bivariatebijection}. If $B^T$ is also partitioned into $k\times k$ blocks, then its blocks are $(B^T)_{ij} = (B_{ji})^T \in \F^{n\times m}$, which implies that its associated $n\times m$ bivariate matrix polynomial is $F(y,x)^T$. Therefore, \eqref{eq:bezout} implies $(\mathbb{DL}(P,v))^T \mapsto \mathcal{B}_{P,v}(y,x,z)^T = \mathcal{B}_{P,v}(x,y,z)^T$. Moreover, from \eqref{eq:bezout},
\[
\mathcal{B}_{P,v} (x,y,z)^T = \frac{P(y)^T (x-z)v(x)-P(x)^T (y-z)v(y)}{x-y} = \mathcal{B}_{P^T,v} (x,y,z).
\]
Thus, $(\mathbb{DL}(P,v))^T$ and $\mathbb{DL}(P^T,v)$ are mapped into the same image via the bijection \eqref{eq:bivariatebijection} and, therefore, $(\mathbb{DL}(P,v))^T = \mathbb{DL}(P^T,v)$.
\end{proof}

We conclude this section by using Theorem \ref{blockevaluation} to obtain for each pencil in the $\mathbb{DL}(P)$ vector space another pencil strictly equivalent to it with a remarkable simple block-sparsity pattern. This indicates that Theorem \ref{blockevaluation} is a relevant result that reveals some intrinsic structure of the pencils in $\mathbb{DL} (P)$.

\begin{remark}
As a side product of Theorem \ref{blockevaluation}, we obtain a new block-arrowhead pencil which is strictly equivalent to $\mathbb{DL}(P,v)$. Indeed, it is obvious that any pencil $A(z) = A_1 z + A_0$ can be written as $A(z) = A(\mu_0) + (z-\mu_0) A_1$, for any arbitrary choice of $\mu_0$. Assume, for maximal sparsity, that $v(z)$ in Theorem \ref{blockevaluation} has finite distinct roots $\mu_1,\mu_2,\dots,\mu_{k-1}$ all with multiplicity one and let $\mu_0 \ne \mu_i$ for $i=1,\ldots, k-1$. Consider $A(z) = A_1 z + A_0 = (W^T \otimes I_m) L(z) (W \otimes I_n)$, where recall that $W$ is nonsingular. Then, by Theorem \ref{blockevaluation},
\[  A(\mu_0) = v(\mu_0) P(\mu_0) \oplus \bigoplus_{i=1}^{k-1} (\mu_i-\mu_0) v'(\mu_i) P(\mu_i) \]
is block diagonal when viewed as a $k\times k$ matrix with blocks of size $m\times n$. Next, we prove that all the blocks of $A_1$ outside the main block-diagonal, the first block-column and the first block-row are identically zero and, moreover, we find explicit expressions for the possibly non-zero blocks. For this purpose, observe that according to \eqref{eq:bivariatebijection} and \eqref{eq:bezout}, the $(j,i)$th block of $A_1$ is prescribed by the formula
\[   \frac{v(y)P(x)-P(y)v(x)}{x-y}  \]
followed by evaluation at $x=\mu_i$, $y=\mu_j$, for $0\leq i,j \leq k-1$. Suppose first $i=j \neq 0$, then by L'Hospital rule we obtain the block $-v'(\mu_i) P(\mu_i)$. If $0 \neq i \neq j \neq 0$, we obtain the zero block. If $i=j=0$, we get $v(\mu_0) P'(\mu_0) - P(\mu_0) v'(\mu_0)$. If $i\ne 0 = j$, we obtain
\[  \frac{v(\mu_0) P(\mu_i)}{\mu_i- \mu_0} ,  \]
and if $i = 0 \ne j$, we obtain the same expression replacing $i$ by $j$.

Applying the results above to the case $k=3$, with $v(z) = z^2 -1$, thus $\mu_1=1$ and $\mu_2=-1$, and taking $\mu_0 = 0$, we get that the following pencil
\[   A(z) = \begin{bmatrix}
-P(0) & 0 & 0\\
0 & 2P(1) & 0\\
0 & 0 & 2P(-1)
\end{bmatrix} + z \begin{bmatrix}
-P'(0) & -P(1) & P(-1) \\
 -P(1) & -2P(1) & 0\\
 P(-1) &0  & 2P(-1)
\end{bmatrix}     \]
is strictly equivalent to $\mathbb{DL}(P,v)$ and, therefore, has the same elementary divisors and minimal indices as $\mathbb{DL}(P,v)$. Hence, $A(z)$ is a linearization of $P(z)$ if $P(z)$ is regular and the eigenvalue exclusion theorem holds, or it is a pencil to which the analysis of this paper applies if $P(z)$ is singular and the eigenvalue exclusion theorem holds. The computational cost of constructing this pencil is however nonzero, since it requires to evaluate $P(z)$ and $P'(z)$ at some points.
\end{remark}

 \section{The elementary divisors of $\mathbb{DL}(P,v)$} \label{sec-elemdivisors}
The goal of this section is to prove item 3 of Theorem \ref{thm:main}, i.e., to show that $\lambda \in \F \cup \{\infty\}$ is an eigenvalue of $P(z)$ if and only if $\lambda$ is an eigenvalue of $L(z)=\mathbb{DL}(P,v)$ and that its partial multiplicities as an eigenvalue of $P(z)$ and as an eigenvalue of $L(z)$ are the same, under the assumptions of Theorem \ref{thm:main}. Thanks to Theorem \ref{thm:infinitoacasa}, we can assume in the arguments of this section that $P(z)$ has no infinite eigenvalues and that $v(z)$ has no infinite roots. The main result in this section is Theorem \ref{thm:partmult}. Its proof requires two preliminary results: Proposition \ref{maximaltolambdaindependent} about the relationship between the root polynomials of $P(z)$ and $L(z)$ and Lemma \ref{embeh} about the orders of the elements of a $\lambda$-independent set of root polynomials of an arbitrary matrix polynomial.

\begin{proposition}\label{maximaltolambdaindependent} Let $P(z) \in \mathbb{F}[z]^{m\times n}$ be a matrix polynomial of grade $k\geq 2$ and $v(z)$ be a scalar ansatz polynomial of grade $k-1$. Suppose that the sets of eigenvalues of $P(z)$ and of roots of $v(z)$ are disjoint, that P(z) has no infinite eigenvalues and that v(z) has no infinite roots. Let $L(z) = \mathbb{DL}(P, v) \in \mathbb{F}[z]^{km \times kn}$ and $V(z) \in \mathbb{F}[z]^{k\times 1}$ be the Vandermonde vector of degree $k-1$ defined in Lemma \ref{important}.
Then:
\begin{enumerate}
\item If $r(z) \in \mathbb{F}[z]^{n\times 1}$ is a root polynomial for $P(z)$ of order $\ell$ at $\lambda$, then $\rho(z)= V(z) \otimes r(z) \in \mathbb{F}[z]^{kn\times 1}$ is a root polynomial for $L(z)$ of order $\ell$ at $\lambda$. In particular, if $\lambda$ is an eigenvalue of $P(z)$, then $\lambda$ is an eigenvalue of $L(z)$.
\item If $\{ r_i(z) \}_{i=1}^t$ is a maximal set of root polynomials at $\lambda$ for $P(z)$, then $\{ \rho_i(z) \}_{i=1}^t$, with $\rho_i(z) = V(z) \otimes r_i(z)$, is a $\lambda$-independent set of root polynomials for $L(z)$.
\end{enumerate}
\end{proposition}
\begin{proof} Item 1.
By assumption, $P(z)r(z)=(z-\lambda)^\ell s(z)$ with $s(\lambda) \neq 0$ and, moreover, $r(\lambda) \not \in \ker_\lambda P(z) = \mathrm{span} \, M(\lambda)$, where $M(z) \in \F [z]^{n\times p}$ is a right minimal basis for $P(z)$.

Observe first that \eqref{eq:bivariatebijection} and \eqref{eq:bezout} imply
$$ L(z) \rho(z) \mapsto \left. \frac{P(y)(x-z)v(x)-P(x)(y-z)v(y)}{x-y} \right|_{x=z} r(z) = v(y) P(z) r(z) = v(y) (z-\lambda)^\ell s(z).$$
If $v(y) = \omega^T \, V(y) = V(y)^T \, \omega$, with $\omega \in \F^{k\times 1}$, the previous equation implies in turn that
\[
L(z) \rho(z) = (z-\lambda)^\ell \, (\omega \otimes s(z)).
\]
Since the ansatz polynomial $v(y)$ is not zero, we get that $\omega \ne 0$ and that $\omega \otimes s(\lambda) \ne 0$. Therefore, it only remains to prove that $\rho(\lambda) \not \in \ker_\lambda L(z)$. 
 For this purpose, note first that from Theorem \ref{thm:Fisminimal} we get $\ker_\lambda L(z) = \spann F(\lambda)$. Moreover, the fact that $r(z)$ is a root polynomial at $\lambda$ for $P(z)$ implies that $\lambda$ is a finite eigenvalue of $P(z)$, since root polynomials at $\lambda$ for $P(z)$ exist if and only if $\lambda$ is a finite eigenvalue of $P(z)$ \cite[Proposition 2.12]{DopicoNoferini}. Thus, $v(\lambda) \ne 0$ by assumption and the matrix $S(\lambda)$ in \eqref{eq.factorizationFz} is well-defined and is nonsingular. This implies that $F(\lambda) = G(\lambda) \, S(\lambda)$, where $G(z) = \begin{bmatrix} C & V(z) \otimes M(z) \end{bmatrix}$ is the matrix polynomial defined in Proposition \ref{thetrick}. Therefore,
\begin{equation} \label{eq.equalspanFG}
\ker_\lambda L(z) = \spann F(\lambda) = \spann G (\lambda).
\end{equation}
Next, we proceed by contradiction and assume that $\rho(\lambda) = V(\lambda) \otimes r (\lambda) \in \ker_\lambda L(z) = \spann G (\lambda)$. Then there exists a constant vector $d \in \F^{kp \times 1}$ such that \begin{equation}\label{eq:abcvr}
\begin{bmatrix} C & V(\lambda) \otimes M(\lambda) \end{bmatrix} d = V(\lambda) \otimes r(\lambda).
\end{equation}
If we partition $d$ conformably to $\begin{bmatrix} C & V(\lambda) \otimes M(\lambda) \end{bmatrix}$ as $d = \left[ \begin{smallmatrix} d_C \\d_M
\end{smallmatrix}\right]$, then \eqref{eq:abcvr} is equivalent to
\[
C \, d_C + V(\lambda) \otimes (M(\lambda)\, d_M) = V(\lambda) \otimes r(\lambda) \Longleftrightarrow \begin{bmatrix} C & V(\lambda) \otimes I_n \end{bmatrix} \begin{bmatrix}
d_C   \\
M(\lambda)\, d_M - r(\lambda)
\end{bmatrix} = 0.
\]
The matrix $\begin{bmatrix} C & V(\lambda) \otimes I_n \end{bmatrix}$ has full column rank, as a consequence of the definition of the matrix $C$ in Proposition \ref{thetrick}, of $\lambda \ne \mu_i$ for $i=1,\ldots , s$, and of Lemma \ref{confvanderparttwo} (with $M_{s+1} (z) = I_n$, $\ell_{s+1} = 1$ and $\alpha_{s+1} = \lambda$). Therefore, $d_C = 0$ and
$M(\lambda)\, d_M - r(\lambda) = 0$, which implies $r(\lambda) \in \ker_\lambda P(z) = \mathrm{span} \, M(\lambda)$, contradicting the assumption.

%In \eqref{eq:abcvr}, we use the notation of Lemma \ref{important} to denote by %$V(\lambda)$ the Vandermonde vector at $\lambda$; moreover, we also use the %notation of Lemmata \ref{confvanderpartone} and \ref{confvanderparttwo}, so that
%$A=\begin{bmatrix}
%A_1 & A_2 & \dots A_m & V(\lambda) \otimes I_n
%\end{bmatrix}$ and $A_i  = \begin{bmatrix}
%V^{(\ell_i-1)} (\mu_i) & \dots & V'(\mu_i) & V(\mu_i)
%\end{bmatrix} \otimes I_n$, where $\mu_i$ is the $i$th root of $v(x)$, having %multiplicity $\ell_i$. Finally, still following the same notation, %$B=\bigoplus_{i=1}^m B_i \oplus M(\lambda)$, and each diagonal block has full %column rank.

%Note that, by construction , $A$ is a square invertible confluent Vandermonde %matrix. Its inverse $A^{-1}$ is the matrix whose rows display the coefficient of %the corresponding Hermite interpolating polynomial at $\mu_i$ ($i=1,\dots,m$), %ith multiplicity $\ell_i$, and $\lambda$. In particular, $A^{-1} V(\lambda)=e_n$. %Thus, $AB c = V(\lambda) \otimes r(\lambda)$ implies $B c = e_n \otimes %r(\lambda)$. Partitioning
%$$ c = \begin{bmatrix}
%c_1\\
%\vdots\\
%c_m\\
%c_\lambda
%\end{bmatrix}$$ coherently with the direct sum structure of $B$, this in turn %implies
%$$ \begin{cases}
%B_1 c_1 = 0 \Rightarrow c_1 = 0;\\
%\vdots
%B_m c_m = 0 \Rightarrow c_m = 0;\\
%M(\lambda) c_\lambda = r(\lambda).
%\end{cases}.$$
%But then, $r(\lambda) \in \ker_\lambda M(\lambda)$, contradicting the assumption.

Item 2. By item 1, it is clear that $\rho_i(z)$ are root polynomials at $\lambda$ for $L(z)$. It remains to check the claim of $\lambda$-independence. Let
$$ R(z) = \begin{bmatrix}
r_1(z) & \cdots & r_t(z)
\end{bmatrix}.$$
Then, taking into account \eqref{eq.equalspanFG} and $G(\lambda) = \begin{bmatrix} C & V(\lambda) \otimes M(\lambda) \end{bmatrix}$, the goal is to prove that the matrix
$$\begin{bmatrix}
C & V(\lambda) \otimes M(\lambda) & V(\lambda) \otimes R(\lambda)
\end{bmatrix} =
\begin{bmatrix}
C & V(\lambda) \otimes \begin{bmatrix} M(\lambda) &  R(\lambda) \end{bmatrix}
\end{bmatrix} \,
$$
has full column rank. But this follows from Lemma \ref{confvanderparttwo}, with
$M_{s+1} (z) = \begin{bmatrix} M(z) &  R(z) \end{bmatrix}$, $\ell_{s+1} = 1$ and $\alpha_{s+1} = \lambda$, and the fact that $ \begin{bmatrix} M(\lambda) &  R(\lambda) \end{bmatrix}$ has full column rank, since $\{ r_i(z) \}_{i=1}^t$ is a maximal set of root polynomials at $\lambda$ for $P(z)$.
\end{proof}

\begin{lemma}\label{embeh}
Let $\{ \ell_i \}_{i=1}^t$ be the partial multiplicities at $\lambda \in \F$ of a matrix polynomial $Q(z)\in \F [z]^{m\times n}$, and let $\{ r_i(z) \}_{i=1}^c$ be a $\lambda$-independent set of root polynomials at $\lambda$ for $Q(z)$ of orders $\{ \kappa_i \}_{i=1}^c$. Then, $c \leq t$ and if both $\{\ell_i\}_{i=1}^t$ and $\{\kappa_i\}_{i=1}^c$ are listed in  non-increasing order, it holds $\kappa_i \leq \ell_i$ for $i=1,\dots,c$.
\end{lemma}
\begin{proof} The definition and properties of $\lambda$-independent and maximal sets of root polynomials at $\lambda$ for $Q(z)$ imply immediately that $c \leq t$. Let $M(z)$ be a minimal basis for $Q(z)$ and let $v_1,\dots,v_{t-c}$ be vectors in $\F^n$ that complete the columns of $\begin{bmatrix}
M(\lambda) & r_1(\lambda) & \cdots & r_c(\lambda)
\end{bmatrix}$ to a basis for $\ker Q(\lambda) \subseteq \F^n$. Then, it is clear that $\{ r_i(z) \}_{i=1}^c \cup \{ v_i \}_{i=1}^{t-c}$ is a complete set of root polynomials at $\lambda$ for $Q(z)$, since $Q(\lambda) v_i = 0$ and $v_i \notin \spann M(\lambda)$. Denote their orders, listed in a non-increasing manner, by $\{\hat \ell_i\}_{i=1}^t$. By Theorem \ref{whenismaximal}, we have $\hat \ell_i \leq \ell_i$ for $i=1,\dots,t$. On the other hand, since the integers $\{\kappa_i\}_{i=1}^c$ are by construction a subsequence of $\{\hat \ell_i \}_{i=1}^t$, the inequalities $\kappa_i \leq \hat \ell_i \leq \ell_i$ must hold for $i=1,\ldots , c$.
\end{proof}

Next, we prove item 3 of Theorem \ref{thm:main}. The proof is based on Theorem \ref{thm:Parts1and2main}, the index sum theorem \cite[Theorem 6.5]{DDM14}, Proposition \ref{maximaltolambdaindependent} and Lemma \ref{embeh}.
\begin{theorem}\label{thm:partmult} {\rm (Item 3 of Theorem \ref{thm:main})} Let $P(z) \in \F[z]^{m \times n}$ be a matrix polynomial of grade $k\geq 2$ and let $v(z) \in \F[z]$ be a scalar polynomial of grade $k-1$. Suppose that the sets of eigenvalues of $P(z)$ and of roots of $v(z)$, each including possibly infinite eigenvalues or roots, are disjoint. Denote by $L(z)\in \F [z]^{km \times kn}$ the pencil in $\mathbb{DL}(P)$ with ansatz polynomial $v(z)$, i.e., $L(z) = \mathbb{DL} (P,v)$. Then, $\lambda \in \F \cup \{\infty\}$ is an eigenvalue of $P(z)$ if and only if $\lambda$ is an eigenvalue of $L(z)$ and, moreover, the partial multiplicities of $\lambda$ as an eigenvalue of $P(z)$ and as an eigenvalue of $L(z)$ coincide.
\end{theorem}
\begin{proof} Similarly to the proof of Theorem \ref{thm:Parts1and2main}, taking into account Theorems \ref{thm:infinitoacasa} and \ref{thm:eigenstrucutreafterMobius}, we can assume without loss of generality that $P(z)$ has no infinite eigenvalues and that $v(z)$ has no infinite roots, for if not we can just apply an appropriate M\"{o}bius transformation to $P(z)$ and $v(z)$. The proof follows an argument by exhaustion. By Theorem \ref{thm:Parts1and2main}, and with the notation used there, we obtain that
$$
\mathrm{grade}(L) \, \rank (L) = kn - kp = k (n-p) = \mathrm{grade}(P) \, \rank (P),
$$
since $\mathrm{grade}(L) = 1$. Moreover, from Theorem \ref{thm:Parts1and2main}, we know that the sum of the right and the left minimal indices of $L(z)$ and $P(z)$ are the same. Hence, by the index sum theorem \cite[Theorem 6.5]{DDM14},
\begin{equation}\label{eq:quantipiccioni}
 \sum_{\lambda \in \Lambda(P)} \sum_{i=1}^{g(\lambda,P)} \ell_i(\lambda,P) = \sum_{\lambda \in \Lambda(L)} \sum_{i=1}^{g(\lambda,L)} \ell_i(\lambda,L),
\end{equation}
where $\Lambda(P)$ denotes the set of distinct eigenvalue of $P(z)$ (and analogously for $\Lambda(L)$), $g(\lambda,P)$ is the geometric multiplicity of $\lambda$ as an eigenvalue of $P$ (and similarly for $g(\lambda,L)$), and $\ell_i(\lambda,P)$ is the $i$-th largest partial multiplicity of $\lambda$ as an eigenvalue of $P$ (and ditto for $\ell_i(\lambda,L)$). Recall that the partial multiplicities have been defined to be positive integers and therefore $g(\lambda,P)$ is equal to the number of partial multiplicities of $\lambda$ as an eigenvalue of $P(z)$.

Observe that the following results hold:
\begin{itemize}
\item[(a)] $\Lambda(P) \subseteq \Lambda(L)$, by item 1 in Proposition \ref{maximaltolambdaindependent},
\item[(b)] $g(\lambda,P) \leq g(\lambda,L)$ for all $\lambda \in \Lambda(P)$, by item 2 in Proposition \ref{maximaltolambdaindependent} and the properties of maximal sets of root polynomials, and
\item[(c)] $\ell_i(\lambda,P) \leq \ell_i(\lambda,L)$ for all $\lambda \in \Lambda(P)$ and for $i = 1, \ldots, g(\lambda,P)$, by Lemma \ref{embeh} and Proposition \ref{maximaltolambdaindependent}.
\end{itemize}
As a consequence, we have that at each $\lambda \in \Lambda(P)$ it must hold
\begin{equation}\label{eq:consequence1}
 \sum_{i=1}^{g(\lambda,P)} \ell_i(\lambda,P) \leq  \sum_{i=1}^{g(\lambda,P)} \ell_i(\lambda,L) \leq \sum_{i=1}^{g(\lambda,L)} \ell_i(\lambda,L),
\end{equation}
where the left inequality is strict if and only if $\ell_i(\lambda,P) < \ell_i(\lambda,L)$ for at least one $i = 1, \ldots , g(\lambda, P)$ and the right inequality is strict if and only if $g(\lambda,P) < g(\lambda,L)$.

We prove first that $\Lambda (P) = \Lambda (L)$ by contradiction. Assume that $\Lambda (P) \subset \Lambda (L)$, i.e., the inclusion in (a) above is strict. Then, from \eqref{eq:quantipiccioni},
\begin{align*}
\sum_{\lambda \in \Lambda(P)} \sum_{i=1}^{g(\lambda,P)} \ell_i(\lambda,P) & = \sum_{\lambda \in \Lambda(P)} \sum_{i=1}^{g(\lambda,L)} \ell_i(\lambda,L) +
\sum_{\lambda \in \Lambda(L)\setminus \Lambda(P)} \sum_{i=1}^{g(\lambda,L)} \ell_i(\lambda,L) \\
& > \sum_{\lambda \in \Lambda(P)} \sum_{i=1}^{g(\lambda,L)} \ell_i(\lambda,L) \\
& \geq \sum_{\lambda \in \Lambda(P)} \sum_{i=1}^{g(\lambda,P)} \ell_i(\lambda,P),
\end{align*}
where the first strict inequality follows from the strict inclusion $\Lambda (P) \subset \Lambda (L)$ and the last inequality follows from \eqref{eq:consequence1}. %Thus, $\sum_{\lambda \in \Lambda(P)} \sum_{i=1}^{g(\lambda,P)} \ell_i(\lambda,P) > \sum_{\lambda \in \Lambda(P)} \sum_{i=1}^{g(\lambda,P)} \ell_i(\lambda,P)$, which is not possible. 
This is a manifest contradiction. Therefore, $\Lambda (P) = \Lambda (L)$.

Secondly, we prove that $g(\lambda,P) = g(\lambda,L)$ for all $\lambda \in \Lambda (P) = \Lambda (L)$. We proceed again by contradiction. If $g(\lambda_0,P) < g(\lambda_0,L)$ for some $\lambda_0 \in \Lambda (P) = \Lambda (L)$, then the second inequality in \eqref{eq:consequence1} is strict for $\lambda_0$ and $\sum_{i=1}^{g(\lambda_0,P)} \ell_i(\lambda_0,P) < \sum_{i=1}^{g(\lambda_0,L)} \ell_i(\lambda_0,L)$. Combining this with \eqref{eq:quantipiccioni}, \eqref{eq:consequence1},  and $\Lambda (P) = \Lambda (L)$,  we get
$$
\sum_{\lambda \in \Lambda(P)} \sum_{i=1}^{g(\lambda,L)} \ell_i(\lambda,L)> \sum_{\lambda \in \Lambda(P)} \sum_{i=1}^{g(\lambda,P)} \ell_i(\lambda,P) = \sum_{\lambda \in \Lambda(P)} \sum_{i=1}^{g(\lambda,L)} \ell_i(\lambda,L),
$$
which is not possible. Then, $g(\lambda,P) = g(\lambda,L)$ for all $\lambda \in \Lambda (P) = \Lambda (L)$.

Finally, we prove $\ell_i(\lambda,P) = \ell_i(\lambda,L)$ for all $\lambda \in \Lambda (P) = \Lambda (L)$ and for all $i = 1,\ldots, g(\lambda,P)= g(\lambda,L)$. We proceed again by contradiction. Assume that there exist at least one $\lambda_0 \in \Lambda (P) = \Lambda (L)$ and at least one $i = 1,\ldots, g(\lambda_0,P)= g(\lambda_0,L)$, such that $\ell_i(\lambda_0,P) < \ell_i(\lambda_0,L)$. Thus, the first inequality in \eqref{eq:consequence1} is strict for $\lambda_0$, i.e., $\sum_{i=1}^{g(\lambda_0,P)} \ell_i(\lambda_0,P) <  \sum_{i=1}^{g(\lambda_0,P)} \ell_i(\lambda_0,L)$. Combining this with \eqref{eq:consequence1}, \eqref{eq:quantipiccioni}, and $g(\lambda,P) = g (\lambda,L)$ for all $\lambda \in \Lambda (P) = \Lambda (L)$, we get
$$
\sum_{\lambda \in \Lambda(P)} \sum_{i=1}^{g(\lambda,P)} \ell_i(\lambda,L)> \sum_{\lambda \in \Lambda(P)} \sum_{i=1}^{g(\lambda,P)} \ell_i(\lambda,P) = \sum_{\lambda \in \Lambda(P)} \sum_{i=1}^{g(\lambda,P)} \ell_i(\lambda,L),
$$
which is not possible. Thus, $\ell_i(\lambda,P) = \ell_i(\lambda,L)$ for all $\lambda \in \Lambda (P) = \Lambda (L)$ and for all $i = 1,\ldots, g(\lambda,P)= g(\lambda,L)$ and the proof is complete.
\end{proof}

We finish this section with Theorem \ref{thm:sec5fin}, which is a corollary of Proposition \ref{maximaltolambdaindependent} and Theorem \ref{thm:partmult}.

\begin{theorem} \label{thm:sec5fin} Let $P(z) \in \mathbb{F}[z]^{m\times n}$ be a matrix polynomial of grade $k\geq 2$ and $v(z)$ be a scalar ansatz polynomial of grade $k-1$. Suppose that the sets of eigenvalues of $P(z)$ and of roots of $v(z)$ are disjoint, that $P(z)$ has no infinite eigenvalues, and that v(z) has no infinite roots. Let $L(z) = \mathbb{DL}(P, v) \in \mathbb{F}[z]^{km \times kn}$ and $V(z) \in \mathbb{F}[z]^{k\times 1}$ be the Vandermonde vector defined in Lemma \ref{important}. Let $\lambda$ be an eigenvalue of $P(z)$ and $\{ r_i(z) \}_{i=1}^t$ be a maximal set of root polynomials at $\lambda$ for $P(z)$. Then, $\{ \rho_i(z) \}_{i=1}^t$, $\rho_i(z) = V(z) \otimes r_i(z)$, is a maximal set of root polynomials at $\lambda$ for $L(z)$.
\end{theorem}

\begin{proof} Proposition \ref{maximaltolambdaindependent} implies that $\{ \rho_i(z) \}_{i=1}^t$ is a $\lambda$-independent set of root polynomials for $L(z)$ with the same orders as $\{ r_i(z) \}_{i=1}^t$. From Theorem \ref{thm:partmult} and the properties of maximal sets of root polynomials, we deduce that the orders of $\{ \rho_i(z) \}_{i=1}^t$ are precisely the partial multiplicites of $\lambda$ as an eigenvalue of $L(z)$. This implies, in particular, that $t = \dim \ker L(\lambda) - \dim \ker_\lambda L(z)$. Therefore, $\{ \rho_i(z) \}_{i=1}^t$ is a complete set of root polynomials at $\lambda$ for $L(z)$. Then, we get that $\{ \rho_i(z) \}_{i=1}^t$ is maximal from Theorem \ref{whenismaximal}.
\end{proof}

% \section{Proof of Theorem \ref{thm:main}}

%We finally have all the ingredients ready to cook up a proof of Theorem %\ref{thm:main}.
% \begin{proof}[Proof of Theorem \ref{thm:main}]
%Suppose first that neither $P(z)$ nor $v(z)I$ have infinite eigenvalues. Then, %for each item:
%\begin{enumerate}
%\item Theorem \ref{thm:abasis} and Theorem \ref{thm:Fisminimal} imply the %statement;
%\item The statement follows from the previous item, by noting that %$\mathbb{DL}(P^T,v)=(\mathbb{DL}(P,v))^T$;
%\item This is a consequence of Theorem \ref{thm:partmult} and Theorem %\ref{whenismaximal}.
%\end{enumerate}
%Finally, once the theorem is established for pairs $(P,v)$ that have no infinite %eigenvalues or roots, the general case follows invoking Theorem %\ref{thm:infinitoacasa} after noting that, since $\F$ is an infinite fields, we %can always construct a M\"{o}bius transforrmation such that neither %$\mathcal{M}_k(P)$ nor $\mathcal{M}_{k-1}(v) I$ have infinite eigenvalues.
% \end{proof}

\section{Recovery of vectors associated with $P(z)$ from their analogues associated with $\mathbb{DL}(P,v)$} \label{sec.recovery}
In this section we describe how to recover (right) minimal bases, eigenvectors, and root polynomials of a possibly singular matrix polynomial $P(z)$ from those of $\mathbb{DL}(P,v)$, under the assumptions of Theorem \ref{thm:main}. We omit the treatment of left minimal bases, eigenvectors, and root polynomials; indeed, it is completely analogous since $\mathbb{DL}(P^T,v)=(\mathbb{DL}(P,v))^T$.

 \subsection{Recovery of a minimal basis for $P(z)$ from one of $\mathbb{DL}(P,v)$}

When studying singular matrix polynomials, often not only the minimal indices are sought, but also a minimal basis. Minimal bases of singular pencils can be computed by postprocessing the output of the staircase algorithm \cite{NVD21,VD-kron-1979}. In this section, we show how a minimal basis for $P(z)$ can be extracted from one of  $L(z)=\mathbb{DL}(P,v)$, under the assumption that the sets of roots of $v(z)$ and of eigenvalues of $P(z)$ are disjoint. It is worth reminding that a minimal basis $F(z)$ of $L(z)$ was built in Theorem \ref{thm:Fisminimal} starting from any minimal basis $M(z)$ of $P(z)$, under the additional assumption%s that $P(z)$ has no infinite eigenvalues and
 that $v(z)$ has no infinite roots. However, there are two obstacles for recovering $M(z)$ from this construction: first $M(z)$ is really ``hidden'' inside $F(z)$ and, second and more important, there may be minimal basis of $L(z)$ with a structure different from the one of $F(z)$. This is easy to see if some of the minimal indices of $P(z)$ are zero. Therefore, we follow in this section a fully different approach.

The vector $\omega \in \F^{k\times 1}$ associated to the scalar polynomial $v(x)$ under the bijection \eqref{eq.omegatov} plays a key role in this section. Recall that $\omega$ is what was called the \emph{ansatz vector} in \cite{MacMMM06}, and it maps under \eqref{eq.omegatov} to what we call, following \cite{NNT17}, the \emph{ansatz polynomial} in this paper. We often use in this section that $v(x) = \omega^T V(x)$, where $V(x)$ is the Vandermonde vector of degree $k-1$ introduced in Lemma \ref{important}.

The following simple lemma will be often used in this section. It follows directly from the definition of $\mathbb{DL}(P)$ pencils and does not require any assumption on $P(z)$ and on $v(z)$.

\begin{lemma} \label{lemm.auxiliaryrecovery} Let $P(z) \in \mathbb{F}[z]^{m\times n}$ be a matrix polynomial of grade $k\geq 2$, $v(z) = \omega^T V(z)$ be a scalar ansatz polynomial of grade $k-1$, where $V(z) \in \F [z]^{k \times 1}$ is the Vandermonde vector of degree $k-1$, and $L(z) = \mathbb{DL}(P, v) \in \mathbb{F}[z]^{km \times kn}$. Then
\begin{enumerate}
\item $(V(z)^T \otimes I_m) \, L (z) = \omega^T \otimes P(z) = P(z) \, (\omega^T \otimes I_n)$;
\item $L(z) \, (V(z) \otimes I_n) = \omega \otimes P(z) = (\omega \otimes I_m) P(z)$;
\item $L(z) \, (V(z) \otimes M(z)) = 0$ for any minimal basis $M(z) \in \F [z]^{n\times p}$ of $P(z)$;
\item $v(z) \, S(z) = (\omega^T \otimes I_n) (V(z) \otimes S(z))$ for any rational matrix (or vector) $S(z) \in \F (z)^{n\times j}$.
\end{enumerate}
The equalities in items 1, 2, and 3 remain valid if the variable $z$ is replaced by any $\lambda \in \F$, and the one in item 4 remains valid if $z$ is replaced by any $\lambda$ that is not a pole of the entries of $S(z)$.
\end{lemma}
\begin{proof}
Item 1 follows from \eqref{eq:bezout}, which implies that $(V(z)^T \otimes I_m) \, L (z) \mapsto \mathcal{B}_{P,v} (x,z,z) = v(x) P(z)$ via the bijection \eqref{eq:bivariatebijection}, while $\omega^T \otimes P(z) \mapsto (\omega^T \otimes P(z)) (V(x) \otimes I_n) = v(x) P(z)$. Analogously, item 2 follows from $L(z) \, (V(z) \otimes I_n) \mapsto \mathcal{B}_{P,v} (z,y,z) = v(y) P(z)$ via the bijection \eqref{eq:bivariatebijection}, while $\omega \otimes P(z) \mapsto (V(y)^T \otimes I_m)(\omega \otimes P(z)) = v(y) P(z)$. Observe that items 1 and 2 follow also from the fact that the set of $\mathbb{DL}(P)$ pencils is the intersection of the sets of $\mathbb{L}_1 (P)$ and $\mathbb{L}_2 (P)$ pencils introduced in \cite{MacMMM06}.

Item 3 was already proved in item 2 of Proposition \ref{thetrick}. It also follows from multiplying on the right the equality in item 2 by $M(z)$ and $P(z) M(z) =0$. Item 4 is an immediate consequence of the properties of the Kronecker product.
\end{proof}

Lemma \ref{lemm.auxiliaryrecovery} allows us to prove the next proposition which, again, does not require any assumption on $P(z)$ and on $v(z)$.

\begin{proposition}\label{prop:recovery}
Let $P(z) \in \mathbb{F}[z]^{m\times n}$ be a matrix polynomial of grade $k\geq 2$, $v(z) = \omega^T V(z)$ be a scalar ansatz polynomial of grade $k-1$, where $V(z) \in \F [z]^{k \times 1}$ is the Vandermonde vector of degree $k-1$, and $L(z) = \mathbb{DL}(P, v) \in \mathbb{F}[z]^{km \times kn}$. Then, the map
$$ \Omega: \ker L(z) \rightarrow \ker P(z), \qquad n(z) \in \F(z)^{kn} \mapsto \Omega(n(z))= (\omega^T  \otimes I_n) \, n(z) \in \F(z)^n $$
is a surjective vector space homomorphism.
\end{proposition}
\begin{proof} Lemma \ref{lemm.auxiliaryrecovery}-1 implies that any $n(z) \in \ker L(z)$ satisfies $0 = (V(z)^T \otimes I_m) \, L (z) \, n(z) = P(z) \, ( (\omega^T \otimes I_n) \, n(z))$. Hence, $(\omega^T \otimes I_n) \, n(z) \in \ker P(z)$. Moreover, $\Omega$ is an $\F(z)$-linear map, and thus a vector space homomorphism from $\ker L(z)$ into $\ker P(z)$. To prove surjectivity, write
an arbitrary $m(z) \in \ker P(z)$ as $m(z)=M(z)\, r(z)$, where $M(z)$ is a minimal basis of $P(z)$ and $r(z)$ is a rational vector. From Lemma \ref{lemm.auxiliaryrecovery}-4, we get $m(z) = v(z) \, M(z) \, (r(z)/v(z)) = (\omega^T \otimes I_n) (V(z) \otimes M(z) \, (r(z)/v(z)))$ and, from Lemma \ref{lemm.auxiliaryrecovery}-3, $L(z) \, (V(z) \otimes M(z) (r(z)/v(z))) = 0$, and the surjectivity of $\Omega$ is proved.
\end{proof}

%To prove surjectivity, our approach is similar to that of the proof of Theorem %\ref{blockevaluation}. Again, let $w(z)=v(z)(z-\mu_0)$, where $\mu_0 \in \F$ does %not belong to the spectrum of $P(z)$ nor to the spectrum of $v(z)$. Define $W$ to %be the (possibly confluent) Vandermonde matrix associated with the roots of %$w(z)$ with their respective multiplicity. Then, we can partition
% $$ W = \begin{bmatrix}
% V(\mu_0) & \hat{W}
% \end{bmatrix}$$
% where $\hat{W} \in \F^{k \times (k-1)}$ is the rectangular, possibly confluent, %Vandermonde matrix associated with the roots of $v(z)$. Observe now that by %construction, there exist a nonzero constant $\beta$ and polynomials %$p_1(z),\dots,p_{k-1}(z)$ such that
%  $$ W^{-1} \begin{bmatrix}
% V(z) & \hat{W}
%\end{bmatrix} = \begin{bmatrix}
%\beta v(z)& 0 & \dots & 0\\
%p_1(z) & 1 & & \\
%\vdots &&\ddots&\\
%p_{k-1}(z)&&&1
%\end{bmatrix}.$$
%On the other hand, $\begin{bmatrix}
%D(z) & C
%\end{bmatrix}$ is a (non-minimal) basis for $\ker L(z)$. Hence, we can write %uniquely $n(z)=D(z) r_1(z) + C r_2(z)$ for some rational $r_1(z),r_2(z)$. By the %observations above, $\beta^{-1} \omega^T = e_1^T W^{-1}$. Thus, %$\Omega(n(z))=(\omega^T \otimes I_n) n(z) = \beta v(z) M(z) r_1(z)$.

%On the other hand, an arbitrary $m(z) \in \ker P(z)$ can be written uniquely as %$m(z)=M(z) r_3(z)$ for some rational $r(z)$. Hence, by picking %$r_1(z)=r_3(z)/(\beta v(z))$ and $r_2(z)$ arbitrarily, we have that $\omega^T %n(z)=m(z)$.
%\end{proof}

From Proposition \ref{prop:recovery} and item 1 of Theorem \ref{thm:Parts1and2main}, we obtain Theorem \ref{thm.recoveryminbases}, which is the main result of this section.

\begin{theorem} (Recovery of minimal bases) \label{thm.recoveryminbases} Let $P(z) \in \F[z]^{m \times n}$ be a matrix polynomial of grade $k \geq 2$, $v(z) =\omega^T V(z)$ be a scalar ansatz polynomial of grade $k-1$, where $V(z) \in \F [z]^{k \times 1}$ is the Vandermonde vector of degree $k-1$, and $L(z) = \mathbb{DL}(P, v) \in \mathbb{F}[z]^{km \times kn}$. Suppose that the sets of eigenvalues of $P(z)$ and of roots of $v(z)$, each including possibly infinite eigenvalues or roots, are disjoint. Let $N(z)$ be a minimal basis of $L(z)$ and define the following matrices:
\begin{enumerate}
\item $\widehat{M}(z) := (\omega^T \otimes I_n) N(z)$;
\item Partition, possibly after having permuted its columns, $\widehat{M}(z)=\begin{bmatrix}
M_d(z) & M_c
\end{bmatrix}$ where no column of $M_d(z)$ is constant in $z$ while $M_c$ is constant;
\item  $\widetilde{M}(z):= M_d(z)$, if $M_c = 0$, and $\widetilde{M}(z):=\begin{bmatrix}
M_d(z) & M_b
\end{bmatrix}$, where $M_b$ is any basis of $\spann M_c$ considered as a subspace over $\F$, otherwise.
\end{enumerate}
Then, $\widetilde{M}(z)$ is a minimal basis of $P(z)$. Thus, $M_d(z)$ is a minimal basis of $P(z)$ if and only if $M_c = 0$.
\end{theorem}
\begin{proof}
Recall that, according to Theorem \ref{thm:Parts1and2main},  $p = \dim \ker P(z)$ if and only if $kp = \dim \ker L(z)$ and that the sum of the (right) minimal indices of $P(z)$ is equal to the sum of those of $L(z)$. Proposition \ref{prop:recovery} implies that the polynomial matrix $\widehat M(z):=(\omega^T \otimes I_n) N(z) \in \F[z]^{n \times kp}$ spans $\ker P(z)$. Moreover, the degrees of the columns of $\widehat M(z)$ are at most equal to the minimal indices of $L(z)$, i.e, to the degrees of the corresponding columns of $N(z)$, which are the minimal indices of $P(z)$ together with $p(k-1)$ minimal indices equal to zero by Theorem \ref{thm:Parts1and2main}. Necessarily, $\widehat M(z)$ is not a basis of $\ker P(z)$ for dimensional reasons, but since its columns are a spanning set, we can construct a polynomial basis (arranged as a polynomial matrix) $B(z)$ of $\ker P(z)$ by selecting $p$ linearly independent columns of $\widehat M(z)$. The fact that any polynomial basis of $\ker P(z)$ must have a sum of the degrees of its vectors greater than or equal to the sum of the minimal indices of $\ker P(z)$ implies that all the positive degrees of the columns of $M_d(z)$ are \emph{equal} to the positive minimal indices of $P(z)$ and that the columns of $M_d(z)$ must be part of the polynomial basis $B(z)$ of $\ker P(z)$, which is in fact minimal due to the value of the sum of the degrees of its vectors. Now, there are two possibilities: (a) $M_d (z)$ has  $p$ columns and (b) $M_d (z)$ has less than $p$ columns. In the case (a) $M_d (z)$ is a minimal basis of $\ker P(z)$ for dimensional reasons and each column of $M_c$ is of the form $M_d(z) x(z)$ for some rational vector $x(z)$, which must be zero, because, if not, then \cite[Main Theorem, Item 4]{For75} would imply that the degree of $M_d(z) x(z)$ is positive.
In the case (b) $B(z) = \begin{bmatrix} M_d (z) & S
\end{bmatrix}$, with $S$ a subset of columns of $M_c$, and each column of $M_c$ is of the form
$\begin{bmatrix} M_d (z) & S
\end{bmatrix} y(z)$ for some rational vector $y(z)$. The entries of $y(z)$ corresponding to $M_d(z)$ must be zero and the ones corresponding to $S$ must be constant, because, if not, then \cite[Main Theorem, Item 4]{For75} would imply again that the degree of $\begin{bmatrix} M_d (z) & S
\end{bmatrix} y(z)$ is positive. Therefore, $\spann_\F S = \spann_\F M_c$, which implies that if $S$ is replaced in $B(z)$ by any other basis $M_b$ of $\spann_\F M_c$, then another minimal basis of $\ker P(z)$ is obtained.
\end{proof}

Theorem \ref{thm.recoveryminbases} can be directly translated into a practical algorithm, where the only step that is not immediate is to construct $M_b$ when $M_c$ is not zero. When $\F \subseteq \mathbb{C}$, which is the most important case in practice, $M_b$ can be reliably computed through a singular value decomposition of $M_c$. (Note that the svd of a matrix over $\mathbb{K}$, a subfield of $\mathbb{C}$, always exists over the closure of $\mathbb{K}$ because an svd can in principle be obtained by solving polynomial equations and performing field operations.) For other fields, one can  apply Gaussian reduction to $M_c$ to identify its pivot columns, which can be taken as $M_b$.

\begin{remark} Note that Proposition \ref{prop:recovery} and Theorem \ref{thm.recoveryminbases} are clearly related to the recovery of the left minimal bases of a {\em singular} matrix polynomial $P(z)$ from those of $\mathbb{L}_1 (P)$ linearizations and to the recovery of the right minimal bases of $P(z)$ from those of $\mathbb{L}_2 (P)$ linearizations presented in \cite[Theorem 5.10]{DDM09}, since in both cases such recoveries are obtained just by the product $(\omega^T \otimes I_n) N(z)$, where $\omega$ is the ansatz vector defining the linearization and $N(z)$ is any of its corresponding minimal bases. Taking into account that the set of $\mathbb{DL}(P)$ pencils is $\mathbb{L}_1 (P) \cap \mathbb{L}_2 (P)$, this relation may seem natural, despite the fact that none of the $\mathbb{DL} (P)$ pencils is a linearization of $P(z)$. However, we emphasize that, in contrast, the recovery of the right minimal bases of $P(z)$ from those of $\mathbb{L}_1 (P)$ linearizations and the recovery of the left minimal bases of $P(z)$ from those of $\mathbb{L}_2 (P)$ linearizations presented in \cite[Theorem 5.10]{DDM09} are not useful at all in the context of $\mathbb{DL}(P)$ pencils for singular $P(z)$.
\end{remark}

Next, we determine in Theorem \ref{thm:kernelOmega} a basis of the kernel of the map $\Omega$ defined in Proposition \ref{prop:recovery}. For simplicity, we use in Theorem \ref{thm:kernelOmega} the additional hypotheses that $P(z)$ has no infinite eigenvalues and that $v(z)$ has no infinite roots. These assumptions can be removed via a M\"{o}bius transformation, but at the cost of making the statement and the proof of the result more cumbersome.

\begin{theorem} \label{thm:kernelOmega} Let $P(z) \in \F[z]^{m \times n}$ be a matrix polynomial of grade $k \geq 2$, $v(z) =\omega^T V(z)$ be a scalar ansatz polynomial of grade $k-1$, where $V(z) \in \F [z]^{k \times 1}$ is the Vandermonde vector of degree $k-1$, $L(z) = \mathbb{DL}(P, v) \in \mathbb{F}[z]^{km \times kn}$, and $\Omega$ be the vector space homomorphism in Proposition \ref{prop:recovery}. Suppose that the sets of eigenvalues of $P(z)$ and of roots of $v(z)$ are disjoint, that $P(z)$ has no infinite eigenvalues and that $v(z)$ has no infinite roots. Let $p = \dim \ker P(z)$ and $C \in \F^{kn \times (k-1)p}$ be the matrix defined in Proposition \ref{thetrick}. Then, $\ker \Omega = \spann C$ and the columns of $C$ form a basis of $\ker \Omega$.
\end{theorem}
\begin{proof} Theorem \ref{thm:Parts1and2main}, Proposition \ref{prop:recovery}, and the rank-nullity theorem imply that $\dim \ker \Omega = p (k-1)$. Moreover, $C$ has full column rank and $\spann C \subseteq \ker L(z)$ by Proposition \ref{thetrick}. Therefore, it only remains to prove that $(\omega^T \otimes I_n) \, C = 0$. For that purpose recall that, according to the definition in Lemma \ref{confvanderpartone}, each submatrix $M_{\mu_i,k,\ell_i,M(z)}$ of $C$ satisfies
\begin{align*}
(\omega^T \otimes I_n) M_{\mu_i,k,\ell_i,M(z)} & =  \left. \begin{bmatrix}
\omega^T  V(z) \otimes M(z) & [\omega^T V(z) \otimes M(z)]' & \cdots & [\omega^T V(z) \otimes M(z)]^{(\ell_i-1)}
\end{bmatrix} \right|_{z=\mu_i} \\
& = \left. \begin{bmatrix}
v(z) M(z) & [v(z) M(z)]' & \cdots & [v(z) M(z)]^{(\ell_i-1)}
\end{bmatrix} \right|_{z=\mu_i} \\
& = 0,
\end{align*}
because $v^{(a)} (\mu_i) = 0$ for $a = 0,1, \ldots , \ell_i -1$.
\end{proof}

Observe that, according to Theorem \ref{thm:kernelOmega}, the minimal indices of the rational subspace $\ker \Omega$ are all equal to zero and that this is the reason why $L(z)$ has $p(k-1)$ minimal indices equal to zero in addition to those of $P(z)$, for any $P(z)$.

We conclude this section with the following result. Again, the assumptions in item 2 of $P(z)$ (resp. $v(z)$) not having infinite eigenvalues (resp. roots) can be removed via a M\"{o}bius transformation.
\begin{corollary}\label{cor:recovery} Let $P(z) \in \mathbb{F}[z]^{m\times n}$ be a matrix polynomial of grade $k\geq 2$, $v(z) = \omega^T V(z)$ be a scalar ansatz polynomial of grade $k-1$, where $V(z) \in \F [z]^{k \times 1}$ is the Vandermonde vector of degree $k-1$, $L(z) = \mathbb{DL}(P, v) \in \mathbb{F}[z]^{km \times kn}$, and $\Omega$ be the vector space homomorphism in Proposition \ref{prop:recovery}. Then,
\begin{enumerate}
  \item $\ker P(z) \cong \ker L(z) / \ker \Omega$;
  \item if, moreover, the sets of eigenvalues of $P(z)$ and of roots of $v(z)$ are disjoint, $P(z)$ has no infinite eigenvalues, $v(z)$ has no infinite roots and $C$ is the matrix defined in Proposition \ref{thetrick}, then $\ker P(z) \cong \ker L(z) / \spann C$.
\end{enumerate}
\end{corollary}
\begin{proof}
Item 1 is an immediate consequence of the first isomorphism theorem \cite[Theorem 6.12]{UA} and the surjectivity of the vector space homomorphism $\Omega$ proved in Proposition \ref{prop:recovery}. Item 2 follows immediately from item 1 and Theorem \ref{thm:kernelOmega}.
\end{proof}

\subsection{Recovery of eigenvectors of $P(z)$ from those of $\mathbb{DL}(P,v)$}
Even though it is not common in the literature to discuss eigenvectors of singular matrix polynomials, it was shown in \cite[\S 2.3]{DopicoNoferini} and in \cite[\S 3]{NVD23} that these can be defined as nonzero elements of the quotient space $\ker P(\lambda)/\ker_\lambda P(z)$ where $\lambda \in \F$ is an eigenvalue of $P(z)$. In other words, given a vector $u$ such that $P(\lambda)u=0$, we say that the equivalence class $[u]=\{u+u' : u' \in \ker_\lambda P(z)\}$ is an eigenvector of $P(z)$ associated with $\lambda$ if $[u] \neq [0]$. This clearly generalizes the notion of an eigenvector of a regular matrix polynomial as a nonzero vector in the null space of the matrix polynomial evaluated at the eigenvalue (indeed, observe that if $P(z)$ is regular then $\ker_\lambda P(z)=\{0\}$). Albeit somewhat abstract, this concept of eigenvectors of singular matrix polynomials has useful applications \cite{KG,LN}. In \cite{KG,LN}, it is also argued that in the (mostly relevant in practice) case $\F \subseteq \C$ then one can make the concept of eigenvectors of a singular matrix polynomial much more concrete. Suppose for example that $\lambda$ is a simple eigenvalue; then this is done by picking, as a representative of the equivalence class of eigenvectors $[u]$, the unique (up to phase) vector $u_0 \in \C^n$ such that (1) $\| u_0 \|_2=1$, (2) $u_0 \in (\ker_\lambda P(z))^\perp$ and (3) $[u_0]=[u]$.

In this sense, a recovery of eigenvectors is possible with a similar approach as for minimal bases. This recovery is based on Proposition \ref{cor:surjection2} and Corollary \ref{cor:recovery2}, which are related to Proposition \ref{prop:recovery} and Corollary \ref{cor:recovery}.

%Indeed, denote as before by $M(z)$ a minimal basis for $\ker P(z)$ and by $R(z)$ %a matrix whose columns are a maximal set of root polynomials at $\lambda$ for %$P(z)$. From the theory of root poynomials \cite{DopicoNoferini} we then know %that for all $\lambda \in \F$ the columns of $\begin{bmatrix}
%M(\lambda) & R(\lambda)
%\end{bmatrix}$ are a basis for $\ker P(\lambda)$, while the columns of %$M(\lambda)$ are a basis for $\ker_\lambda P(z)$. On the other hand, by the %results developed earlier on in this paper, and still under the eigenvalue %exclusion assumptions, for all $\lambda$ such that $v(\lambda)\neq0$ it holds %that $\begin{bmatrix}
%D(\lambda) & C
%\end{bmatrix}$ is a basis for $\ker_\lambda L(z)$ and that $\begin{bmatrix}
%D(\lambda) & C & V(\lambda) \otimes R(\lambda)
%\end{bmatrix}$ is a basis for $\ker P(\lambda)$, implying the following result %whose proof is omitted as it is analogous to Proposition \ref{prop:recovery} %Corollary \ref{cor:recovery}.

\begin{remark} We warn the reader that we are about to commit an abuse of notation, using the same symbol $\Omega$ to denote the map defined in Proposition \ref{prop:recovery}, the two maps defined in Proposition \ref{cor:surjection2} below, and the map defined by \eqref{eq:isomorphism} below. Indeed, the first is an $\F(z)$-linear map defined on $\F(z)^{kn}$ and restricted on $\ker L(z)$; the second and the third are the restrictions of an $\F$-linear map, defined on $\F^{kn}$, on $\ker_\lambda L(z)$ and $\ker L(\lambda)$, respectively; and the fourth is an $\F$-linear map defined on $\ker L(\lambda)/\ker_\lambda L(z)$ (which formally is a quotient $\F$-vector space of equivalence classes). In practice and when represented in the canonical bases, though, all these maps are realized by left-multiplication times $\omega^T \otimes I_n$; moreover, the context always makes clear which one we are referring to. For these reasons, we opt out of an overwhelmingly baroque notation that uses four different symbols.
\end{remark}

\begin{proposition}\label{cor:surjection2} Let $P(z) \in \mathbb{F}[z]^{m\times n}$ be a matrix polynomial of grade $k\geq 2$, $v(z) = \omega^T V(z)$ be a scalar ansatz polynomial of grade $k-1$, where $V(z) \in \F [z]^{k \times 1}$ is the Vandermonde vector of degree $k-1$, and $L(z) = \mathbb{DL}(P, v) \in \mathbb{F}[z]^{km \times kn}$. Fix any $\lambda \in \F$ such that $v(\lambda) \ne 0$. Then,
\begin{enumerate}
\item the map $$\Omega: \ker_\lambda L(z) \rightarrow \ker_\lambda P(z), \qquad u \mapsto (\omega^T \otimes I_n) u$$ is a surjective vector space homomorphism;
\item the map
$$\Omega: \ker L(\lambda) \rightarrow \ker P(\lambda), \qquad u \mapsto (\omega^T \otimes I_n) u$$
is a surjective vector space homomorphism;
\item if, moreover, the sets of eigenvalues of $P(z)$ and of roots of $v(z)$ are disjoint, then the kernel of the map in item 1 is equal to the kernel of the map in item 2, which is denoted by $\ker \Omega$. As a consequence, $\ker_\lambda P(z) \cong \ker_\lambda L(z)/\ker \Omega$ and $\ker P(\lambda) \cong \ker L(\lambda)/\ker \Omega$.
\end{enumerate}
\end{proposition}
\begin{proof}
Item 1. Let us prove first that $\Omega$ indeed maps vectors of $\ker_\lambda L(z)$ to vectors of $\ker_\lambda P(z)$. For this purpose, let $N(z)$ be a minimal basis of $L(z)$ and $M(z)$ be a minimal basis of $P(z)$. Then, Proposition \ref{prop:recovery} implies that $(\omega^T \otimes I_n) N(z) = M(z) R(z)$, where $R(z)$ is a polynomial matrix by \cite[Main Theorem, Item 4]{For75}. Then, $(\omega^T \otimes I_n) N(\lambda) = M(\lambda) R(\lambda)$. By definition, $u \in \ker_\lambda L(z)$ if $u = N(\lambda) s$ for some constant vector $s$. Thus,
$(\omega^T \otimes I_n) u = (\omega^T \otimes I_n) N(\lambda) s = M(\lambda) (R(\lambda) s) \in \ker_\lambda P(z)$. Since $\Omega$ is clearly linear, it only remains to prove its surjectivity. Any $h \in \ker_\lambda P(z)$ can be written as $h = M(\lambda) f$ for some constant vector $f$.  Thus, $h = v(\lambda) M(\lambda) \, (f/ v(\lambda)) = (\omega^T \otimes I_n) (V(\lambda) \otimes  M(\lambda) ) \, (f/ v(\lambda))$, which proves the surjectivity because $(V(\lambda) \otimes  M(\lambda)) \, (f/ v(\lambda)) \in \ker_\lambda L(z)$. This last step follows from Lemma \ref{lemm.auxiliaryrecovery}-3, which implies $V(z) \otimes M(z) = N(z) T(z)$, for some polynomial matrix $T(z)$. Therefore, $(V(\lambda) \otimes  M(\lambda) )\, (f/ v(\lambda)) = N(\lambda)( T(\lambda) \, (f/ v(\lambda))) \in \ker_\lambda L(z)$.

Item 2.  Lemma \ref{lemm.auxiliaryrecovery}-1 with $z$ replaced by $\lambda$ immediately implies that $\Omega$ maps any vector $u \in \ker L(\lambda)$ into a vector in $\ker P(\lambda)$, because $0 = (V(\lambda)^T \otimes I_m) \, L (\lambda) \, u = P(\lambda) \, (\omega^T \otimes I_n) \, u$. To prove surjectivity, write any $h \in \ker P(\lambda)$ as $h = v(\lambda) (h/v(\lambda)) = (\omega^T \otimes I_n) (V(\lambda) \otimes  I_n ) \, (h/ v(\lambda))$. Note that Lemma \ref{lemm.auxiliaryrecovery}-2 with $z$ replaced by $\lambda$ implies that
$L(\lambda) (V(\lambda) \otimes  I_n ) \,(h / v(\lambda)) = (\omega \otimes I_m) P(\lambda) \,(h/ v(\lambda)) = 0$, i.e., $(V(\lambda) \otimes  I_n ) \, (h/ v(\lambda)) \in \ker L(\lambda)$.

Item 3. Theorem \ref{thm:main} implies that $\lambda$ is an eigenvalue of $P(z)$ if and only if $\lambda$ is an eigenvalue of $L(z)$. If $\lambda$ is not an eigenvalue, then $\ker_\lambda P(z) = \ker P(\lambda)$ and $\ker_\lambda L(z) = \ker L(\lambda)$ and there is nothing to prove concerning the equality of the kernels. Thus, we assume that $\lambda$ is an eigenvalue of $P(z)$ and of $L(z)$ and denote the maps in items 1 and 2 by $\Omega_1$ and $\Omega_2$, respectively. Theorem \ref{thm:main} implies that $\dim \ker_\lambda P(z) = p$ if and only if $\dim \ker_\lambda L(z) = pk$ and that $\dim \ker P(\lambda) = p + t$ if and only if $\dim \ker L(\lambda) = pk + t$. The rank-nullity theorem and the surjectivity of $\Omega_1$ and $\Omega_2$ imply that
$\dim \ker \Omega_1 = \dim \ker \Omega_2 = p(k-1)$. Moreover, $\ker_\lambda L(z) \subset \ker L(\lambda)$ implies $\ker \Omega_1 \subseteq \ker \Omega_2$ with strict equality by dimensional reasons. The stated isomorphisms follow from the first isomorphism theorem.
\end{proof}

The third isomorphism theorem \cite[Theorem 6.18]{UA} and item 3 of Proposition \ref{cor:surjection2} immediately imply the following corollary.
\begin{corollary}\label{cor:recovery2}  Let $P(z) \in \F[z]^{m \times n}$ be a matrix polynomial of grade $k \geq 2$, $v(z)$ be a scalar ansatz polynomial of grade $k-1$, and $L(z) = \mathbb{DL}(P, v) \in \mathbb{F}[z]^{km \times kn}$. Suppose that the sets of eigenvalues of $P(z)$ and of roots of $v(z)$ are disjoint  and fix any $\lambda \in \F$  such that $v(\lambda) \ne 0$. Then, $\ker P(\lambda)/\ker_\lambda P(z) \cong \ker L(\lambda)/\ker_\lambda L(z)$.
\end{corollary}

Our proof of the isomorphism in Corollary \ref{cor:recovery2} is abstract in the sense that it reduces to observing that the statement is a special case of deep and basic algebraic results; still, it can be also constructive in the sense of explicitly exhibiting the relevant vector space isomorphism due to items 1 and 2 of Proposition \ref{cor:surjection2}. Namely, the latter is given by the map
\begin{equation}\label{eq:isomorphism}
\Omega :  \ker L(\lambda)/\ker_\lambda L(z) \rightarrow  \ker P(\lambda)/\ker_\lambda P(z), \qquad [u] \mapsto [(\omega^T \otimes I_n) u].
\end{equation}

Finally, observe that, similarly to Theorem \ref{thm:kernelOmega}, one can prove that the kernels of the vector space homomorphisms in items 1 and 2 of Proposition \ref{cor:surjection2} are both equal to $\spann C$.

\subsection{Recovery of root polynomials of $P(z)$ from those of $\mathbb{DL}(P,v)$}
Finally, in this section we prove that the same approach used in the previous subsections for right minimal bases and right eigenvectors (that we can informally describe as ``left multiply times $\omega^T \otimes I_n$'') allows us also to recover right root polynomials. For a pencil, maximal sets of root polynomials can be computed starting from the staircase form of the pencil \cite{NVD21}. Therefore, the result in this subsection can be used in practice to compute maximal sets of root polynomials for $P(z)$. For brevity, we focus on root polynomials at finite eigenvalues.

\begin{theorem} \label{thm.recoveryrootpolys} Let $P(z) \in \F[z]^{m \times n}$ be a matrix polynomial of grade $k \geq 2$, $v(z) =\omega^T V(z)$ be a scalar ansatz polynomial of grade $k-1$, where $V(z) \in \F [z]^{k \times 1}$ is the Vandermonde vector of degree $k-1$, and $L(z) = \mathbb{DL}(P, v) \in \mathbb{F}[z]^{km \times kn}$. Suppose that the sets of eigenvalues of $P(z)$ and of roots of $v(z)$, each including possibly infinite eigenvalues or roots, are disjoint.
Let $\{ \rho_i(z) \}_{i=1}^t$ be a maximal set of root polynomials at $\lambda \in \F$ for $L(z)$ and define $r_i(z) = (\omega^T \otimes I_n) \rho_i(z)$ for $i=1,\dots,t$. Then, $\{ r_i (z) \}_{i=1}^t$ is a maximal set of root polynomials at $\lambda$ for $P(z)$.
\end{theorem}
\begin{proof} Let $p = \dim \ker P(z)$. Then, $pk = \dim \ker L(z)$ by Theorem \ref{thm:main}. Let $N(z) \in \F[\lambda]^{nk \times pk}$ be a minimal basis of $L(z)$ and let $M(z) \in \F[\lambda]^{n \times p}$ be a minimal basis of $P(z)$ obtained by applying Theorem \ref{thm.recoveryminbases} to $N(z)$ with $M_b$ a maximal set of linearly independent columns of $M_c$ (in case $M_c \ne 0$). Thus, $M(\lambda)$ is formed by a subset of the columns of $\widehat{M} (z) = (\omega^T \otimes I_n) N(z)$. By definition, the columns of $N(\lambda)$ form a basis of $\ker_\lambda L(z)$ and the columns of
$$
N_{ext} (\lambda) = \begin{bmatrix}
N(\lambda) & \rho_1 (\lambda) & \cdots & \rho_t (\lambda)
\end{bmatrix} \in \F^{nk \times (pk+t)}
$$
form a basis of $\ker L(\lambda)$. From Proposition \ref{cor:surjection2}, the columns of $\widehat{M} (\lambda) = (\omega^T \otimes I_n) N(\lambda)$ span $\ker_\lambda P(z)$ and the ones of $(\omega^T \otimes I_n) N_{ext} (\lambda)$ span $\ker P(\lambda)$. But those columns of $\widehat{M} (\lambda)$ corresponding to $M(\lambda)$ also span $\ker_\lambda P(z)$, which implies that the columns of
$$
\begin{bmatrix}
M(\lambda) & r_1 (\lambda) & \cdots & r_t (\lambda)
\end{bmatrix} \in \F^{n \times (p+t)}
$$
span $\ker P(\lambda)$. Theorem \ref{thm:main} guarantees that $\dim \ker P(\lambda) = p + t$. Therefore, the columns of $\begin{bmatrix}
M(\lambda) & r_1 (\lambda) & \cdots & r_t (\lambda)
\end{bmatrix}$ form a basis of $\ker P(\lambda)$, which immediately implies that $r_i (\lambda) \notin \ker_\lambda P(z)$ for $i = 1, \ldots , t$. Therefore, it only remains to prove that $r_1 (z) ,\ldots , r_t (z)$ are root polynomials at $\lambda$ for $P(z)$ with the same orders as $\rho_1 (z) ,\ldots , \rho_t (z)$, since by Theorem \ref{thm:main} these orders are the partial multiplicities of $\lambda$ as an eigenvalue of $P(z)$ and Theorem \ref{whenismaximal} would imply the maximality of $\{ r_i (z) \}_{i=1}^t$. For this purpose, note that from the equations $L(z) \rho_i(z) = (z-\lambda)^{\ell_i} s_i(z)$ with $s_i(\lambda) \neq 0$, corresponding to the fact that $\{ \rho_i(z) \}_{i=1}^t$ are root polynomials at $\lambda$ for $L(z)$ of orders $\{\ell_i \}_{i=1}^t$, and from Lemma \ref{lemm.auxiliaryrecovery}-1, we obtain
$$(z-\lambda)^{\ell_i} (V(z)^T \otimes I_m) s_i(z) = (V(z)^T \otimes I_m) L(z) \rho_i(z)= P(z) (\omega^T \otimes I_n) \rho_i(z)$$ and $(V(\lambda)^T \otimes I_m) s_i(\lambda) \ne 0$. Thus, $r_i(z) = (\omega^T \otimes I_n) \rho_i(z)$ is a root polynomial at $\lambda$ for $P(z)$ of order $\ell_i$, for $i=1, \ldots , t$.
\end{proof}

%Denote by $M(z)$ a minimal basis of $\ker P(z)$ and by $$R(z)=\begin{bmatrix}
%r_1(z) & \dots & r_s(z)
%\end{bmatrix}.$$ Note that, by Corollary \ref{cor:recovery2}, the equivalence %classes whose represenatives are the columns of $R(\lambda)$ are a basis for %$\ker P(\lambda)/\ker_\lambda P(z)$ while by Proposition \ref{cor:surjection2} %the columns of $\begin{bmatrix}
%M(\lambda) & R(\lambda)
%\end{bmatrix}$ are a basis for $\ker P(\lambda)$.

%From the equations $L(z) \rho_i(z) = (z-\lambda)^{\ell_i} s(z), \quad s(\lambda) %\neq 0$ we obtain
%$$ (V(z)^T \otimes I_n) L(z) \rho_i(z)= P(z) (\omega^T \otimes I_n) \rho_i(z) = %(x-\lambda)^\ell_i (V(z)^T \otimes I_n) s_i(z).$$
%Suppose that, for some $i$, $s_i(z)=0$. Then $r_i(z) \in \ker P(z)$, which %implies that $\begin{bmatrix}
%M(\lambda) & R(\lambda)
%\end{bmatrix}$ is not full rank: a contradiction since their columns are a basis %of $\ker P(\lambda)$. On other other hand, suppose that for some  $j$ we have %$r_j(\lambda) \in \ker_\lambda P(z)$. Then $R(\lambda)$ is not full rank, another %contradiction since (the equivalence classes represented by) its columns are a %basis for $\ker P(\lambda)/\ker_\lambda P(z)$.

%This shows that $\{ r_i(z) \}_{i=1}^s$ is a complete set of root polynomials at %$\lambda$ for $P(z)$, whose orders are all greater  than or equal to the partial %multiplicities of $\lambda$ in $P(z)$. By Theorem \ref{whenismaximal}, it must %then be a maximal set.

\begin{remark}
The recovery procedure for root polynomials described in Theorem \ref{thm.recoveryrootpolys} is, of course, still valid when $P(z)$ is regular. In that case, root polynomials can be seen as generating functions of Jordan chains \cite{DopicoNoferini}. Thus, this subsection also describes, as a special case, how to recover Jordan chains for a regular $P(z)$ from those of $\mathbb{DL}(P,v)$ under the eigenvalue exclusion assumption. Note that in the regular case, the pencils in $\mathbb{DL} (P)$ satisfying the eigenvalue exclusion condition are linearizations for $P(z)$, then they are also $\mathbb{L}_1 (P)$ and $\mathbb{L}_2 (P)$ linearizations \cite{MacMMM06} and one can also recover maximal sets of root polynomials through a simple block extraction, as it was described in item 1 of \cite[Theorem 8.10]{DopicoNoferini}. In fact, for $P(z)$ regular, Theorem \ref{thm.recoveryrootpolys} follows from item 2 of \cite[Theorem 8.10]{DopicoNoferini}.
\end{remark}

\section{Conclusions} \label{sec-conclusions}

We have extended the eigenvalue exclusion theorem for $\mathbb{DL}(P)$ pencils to the case of a singular matrix polynomial $P(z)$. Even if none of the pencils in $\mathbb{DL}(P)$ is a linearization for a singular $P(z)$, we have proved that, under the eigenvalue exclusion theorem assumptions, it is possible from any pencil $\mathbb{DL}(P,v)$ to fully recover all of the relevant magnitudes of $P(z)$, including eigenvalues and their partial multiplicities, (left and right) eigenvectors, (left and right) root polynomials, (left and right) minimal indices, and (left and right) minimal bases. With the exception of the recovery of the minimal indices and bases, this was already known when $P(z)$ is regular, but the arguments for a singular $P(z)$ are much more involved since they cannot be based on the (no longer true) fact that $\mathbb{DL}(P,v)$ is a strong linearization.

Our analysis also raises a question that goes beyond the study of the space of pencils $\mathbb{DL}(P)$. Does one really need strong linearizations to compute the eigenstructure of a matrix polynomial? Or, is any linear pencil that allows us to recover all the relevant quantities possibly as good as a tool? The case study of $\mathbb{DL}(P)$ definitely suggests that strictly restricting to strong linearizations may be an unnecessarily rigid approach. This is also supported by the new concept of ``strongly minimal linearizations'' introduced in \cite{DQVD-selfconjugate}.

An open problem in this context is to study what happens when the hypotheses of the generalized eigenvalue exclusion theorem are not satisfied. Based on preliminary analyses and also on the results in \cite{DQVD-selfconjugate}, we conjecture that a large part of Theorem \ref{thm:main} will remain valid, but that those eigenvalues of $P(z)$ that are roots of the ansatz polynomial will require of a different analysis and will satisfy different recovery rules for their partial multiplicities.
%\section*{References}

\end{document}